\documentclass[10pt,USletter,twocolumn]{article}
\usepackage{amsmath}
\usepackage{amsfonts}
\usepackage{amssymb}
\usepackage{graphicx}
\usepackage[left=2cm,right=2cm,top=2cm,bottom=2cm]{geometry}

\usepackage{makecell}
\usepackage{forest}
\usepackage{setspace}

\usepackage[electronic]{ifsym}
\usepackage{indentfirst,latexsym}
\usepackage{bm}
\usepackage{mathrsfs}
\usepackage{amsmath,amsfonts,amsthm,amscd,amssymb}
\usepackage{pifont}
\usepackage{url}

\usepackage{longtable}

\usepackage{graphicx}
\usepackage{subfigure}

\usepackage{cite}

\usepackage[all,pdf]{xy}

\usepackage{booktabs}
\usepackage{multirow}
\usepackage{multicol}
\usepackage{stfloats}

\DeclareMathOperator{\dif}{d}  
\DeclareMathOperator{\dist}{\mathrm{dist}}

\renewcommand{\vec}[1]{\mathbf{#1}}

\newcommand{\Real}{\mathbb{R}}

\newcommand{\me} {\mathrm{e}}
\newcommand{\fracode}[2]{\frac{\dif {#1}}{\dif {#2}}}         

\newcommand{\set}[1]{\left\{ #1 \right\}}
\newcommand{\seq}[1]{\left\langle #1 \right\rangle}
\newcommand{\abs}[1]{\left| #1 \right|}
\newcommand{\braket}[2]{ \langle #1 | #2 \rangle}
\newcommand{\norm}[1]{\left\lVert #1 \right\rVert}
\newcommand{\normp}[2]{{\left\lVert #1 \right\rVert}_{#2}}
\newcommand{\trsp}[1]{{#1}^\textsf{T}}

\newcommand{\ES}[3]{\mathbb{#1}^{{#2}\times {#3}}}     

\DeclareMathOperator{\AGM}{AGM}

\newcommand{\scru}[2]{{#1}^{\mathrm{#2}}}
\newcommand{\scrd}[2]{{#1}_{\mathrm{#2}}}

\usepackage{algorithm}         
\usepackage{algorithmicx}
\usepackage{algpseudocode}
\usepackage{float}

\newcommand{\Algr}{\textbf{Algorithm}~}
\newcommand{\Fig}{\textbf{Figure}~}
\newcommand{\Tab}{\textbf{Table}~}

\algnewcommand\algorithmicswitch{\textbf{switch}}
\algnewcommand\algorithmiccase{\textbf{case}}
\algnewcommand\algorithmicdefault{\textbf{default}}

\algdef{SE}[SWITCH]{Switch}{EndSwitch}[1]{\algorithmicswitch\ #1\ \algorithmicdo}{\algorithmicend\ \algorithmicswitch}%
\algdef{SE}[CASE]{Case}{EndCase}[1]{\algorithmiccase\ #1}{\algorithmicend\ \algorithmiccase}%
\algdef{SE}[DEFAULT]{Default}{EndDefault}[1]{\algorithmicdefault\ #1}{\algorithmicend\ \algorithmicdefault}%

\makeatletter
\renewcommand{\ALG@name}{Algorithm}
\newenvironment{breakablealgorithm}
{
	\begin{center}
		\refstepcounter{algorithm}
		\setlength{\baselineskip}{15pt} 
		\renewcommand{\caption}[2][\relax]{
			\hrule height.9pt depth0pt \kern0pt
			{\raggedright\textbf{\ALG@name~\thealgorithm} ##2\par}%
			\ifx\relax##1\relax 
			\addcontentsline{loa}{algorithm}{
				\protect\numberline{\thealgorithm}##2}%
			\else 
			\addcontentsline{loa}{algorithm}{
				\protect\numberline{\thealgorithm}##1}%
			\fi
			\kern2pt\hrule\kern2pt
		}
	}{
		\kern3pt\hrule\relax
	\end{center}
}

\usepackage{caption}
\usepackage{listings}
\newcommand{\cpvar}[1]{\texttt{#1}}

\newcommand{\True}{\texttt{TRUE}}
\newcommand{\False}{\texttt{FALSE}}
\newcommand{\PtrToFun}[1]{\texttt{#1}}
\newcommand{\ProcName}[1]{\textsc{#1}}

\newtheorem{thm}{Theorem}
\newtheorem{cor}[thm]{Corollary}

\author{Hong-Yan Zhang\thanks{Corresponding author, e-mail: hongyan@hainnu.edu.cn}~ 
and Wen-Juan Jiang\\ School of Information Science and Technology, Hainan Normal University, Haikou 571158, China}
\title{Open Source Implementations of Numerical Algorithms for \\ Computing the Complete Elliptic Integral of the First Kind\thanks{The current version is obtained from the earlier version titled with ``Correct and Alternative Numerical Algorithms for the Complete Elliptic Integral of the First Kind in MATLAB and Mathematica" by fixing some bugs and errors.}}
\date{Aug. 12, 2023}
\begin{document}
\maketitle

\begin{abstract}
The complete elliptic integral of the first kind (CEI-1) plays a significant role in mathematics, physics and engineering. There is no simple formula for its computation, thus numerical algorithms are essential for coping with the  practical problems involved. The commercial implementations for the numerical solutions, such as the functions \lstinline|ellipticK| and \lstinline|EllipticK| provided by MATLAB and Mathematica respectively, are based on $\scrd{\mathcal{K}}{cs}(m)$ instead of the usual form $K(k)$ such that $\scrd{\mathcal{K}}{cs}(k^2) =K(k)$ and $m=k^2$. It is necessary to develop open source implementations for the computation of the CEI-1 in order to avoid potential risks of using commercial software and possible limitations due to the unknown factors. In this paper, the infinite series method, arithmetic-geometric mean (AGM) method, Gauss-Chebyshev method and Gauss-Legendre methods are discussed in details with a top-down strategy. The four key algorithms for computing CEI-1 are designed, verified, validated and tested, which can be utilized in R\& D and be reused properly. Numerical results show that our open source implementations based on $K(k)$ are equivalent to the commercial implementation based on $\scrd{\mathcal{K}}{cs}(m)$. The general algorithms for computing orthogonal polynomials developed are significant byproducts in the sense of STEM education and scientific computation. 
\end{abstract}

\tableofcontents

\newpage 

\section{Introduction}  \label{sec-Intro}

The \textit{Complete Elliptic Integral of the first kind} (CEI-1) is an interesting special function in mathematics which also plays a significant role in physics and engineering applications. There are two definitions for the CEI-1: the first one is given   by \cite{TableISP2015,Abramo1965,HYZhang2024CdioCt}
\begin{equation} \label{eq-K-ku}
K(k) = \int^{\pi/2}_0 \frac{\dif u}{\sqrt{1-k^2 \sin^2 u}}, \quad k\in [0,1].
\end{equation}
and the second one is given by
\begin{equation} \label{eq-K-mu}
\scrd{\mathcal{K}}{cs}(m) = \int^{\pi/2}_0 \frac{\dif u}{\sqrt{1-m \sin^2 u}}, \quad m\in [0,1]
\end{equation}
which is implemented by the commercial software, denoted by the subscript \verb|cs| of the function $\scrd{\mathcal{K}}{cs}(\cdot)$,  such as the \lstinline|ellipticK| function in MATLAB and the \lstinline|EllipticK| function in Mathematica. Obviously, the equations \eqref{eq-K-ku} and
\eqref{eq-K-mu} are connected by
\begin{equation} \label{eq-km}
\left\{
\begin{aligned}
& m = k^2 \\
& \scrd{\mathcal{K}}{cs}(m) = K(k)
\end{aligned} 
\right.
\end{equation}
for the arguments $k$ and $m$. 

\begin{figure}[h]
\centering
\subfigure[Output of \lstinline|EllipticK| in Mathematica]{
\boxed{\includegraphics[width=0.45\textwidth]{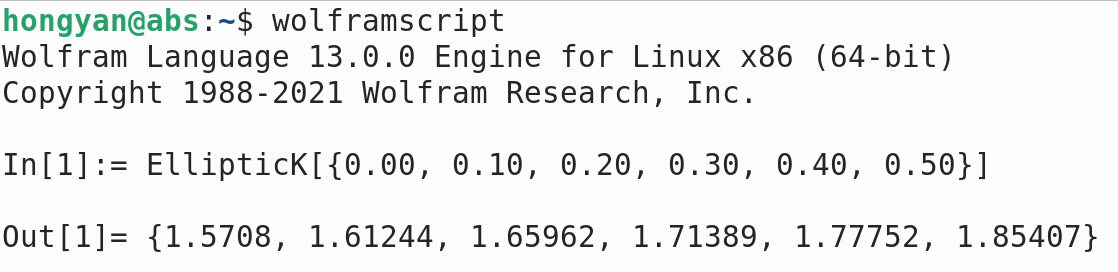}} 
}
\subfigure[Output of \lstinline|ellipticK| in MATLAB]{
\boxed{\includegraphics[width=0.45\textwidth]{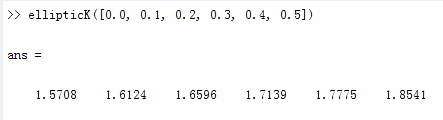}} 
}
\caption{Value of $\scrd{\mathcal{K}}{cs}(m)$ obtained by the MATLAB and Mathematica for $m = [0.00, 0.01, 0.20, 0.30, 0.40, 0.50]$} \label{K-numer-commer}
\end{figure}

For the special function $K(\cdot)$ defined by integral, the computation should be done with proper numerical algorithm. In practical applications, there are different implementations of numerical algorithms for the CEI-1. \Fig \ref{K-numer-commer} illustrates the numerical results obtained by the functions \lstinline|ellipticK| in MATLAB and \lstinline|EllipticK| in Mathematica. However, the commercial implementations are limited by the licenses. Moreover, the application of MATLAB is prohibited for some affiliations due to some potential conflicts involved in national security.

Actually, we can find some fundamental types of methods to compute the CEI-1 properly:
\begin{itemize}
	\item expanding the CEI-1 into an infinite series which converges rapidly;
	\item computing the CEI-1 with the arithmetic-geometric mean (AGM) method which converges rapidly and iterates simply;
	\item converting the CEI-1 into a finite sum with the $2n-1$ order algebraic precision by taking the $n$-th order orthogonal polynomials via Gauss numerical integration approach. 
\end{itemize}	
The Gauss numerical method can be expressed by\cite{Atkinson2009}

\begin{equation}
I(f) = \int^b_a \rho(x)f(x)\dif x 
	= \sum^n_{i=1}w_if(x_i) + \mathcal{E}_n(f)
\end{equation}
where $n$ is a positive integer, $\set{x_i: 0 \le i\le n-1}$ is the set of $n$ roots of an orthogonal polynomial, say $ G_n(x)$, $w_i$ are the weights specified by $x_i$ and $G_n(x)$, and $\mathcal{E}_n(f)$ is the error of optimal approximation which will be omitted according to the requirement of precision in practice. 
	However, the integral $K(k)$ is more complicated since it is a function  defined by an integral. With the help of the concepts of abstract function defined by an integral, we can use the following form	
\begin{equation}
\begin{split}
K(k) &=\int^b_a f(x,k)\dif \mu 
= \int^b_a f(x,k)\rho(x)\dif x \\
&= \sum^n_{i=1} w_i f(x_i, k) + \mathcal{E}_n(f).
\end{split}
\end{equation}
	There are two feasible methods in this way:
	\begin{itemize}
		\item Gauss-Chebyshev method, where $\displaystyle \rho(x) = \frac{1}{\sqrt{1-x^2}}$and $b= -a = 1$;
		\item Gauss-Legendre method, where $\rho(x) = 1$, $a=0$ and $b= \pi/2$.
	\end{itemize}    

\begin{figure}[h]
	\centering
	\includegraphics[width=0.5\textwidth]{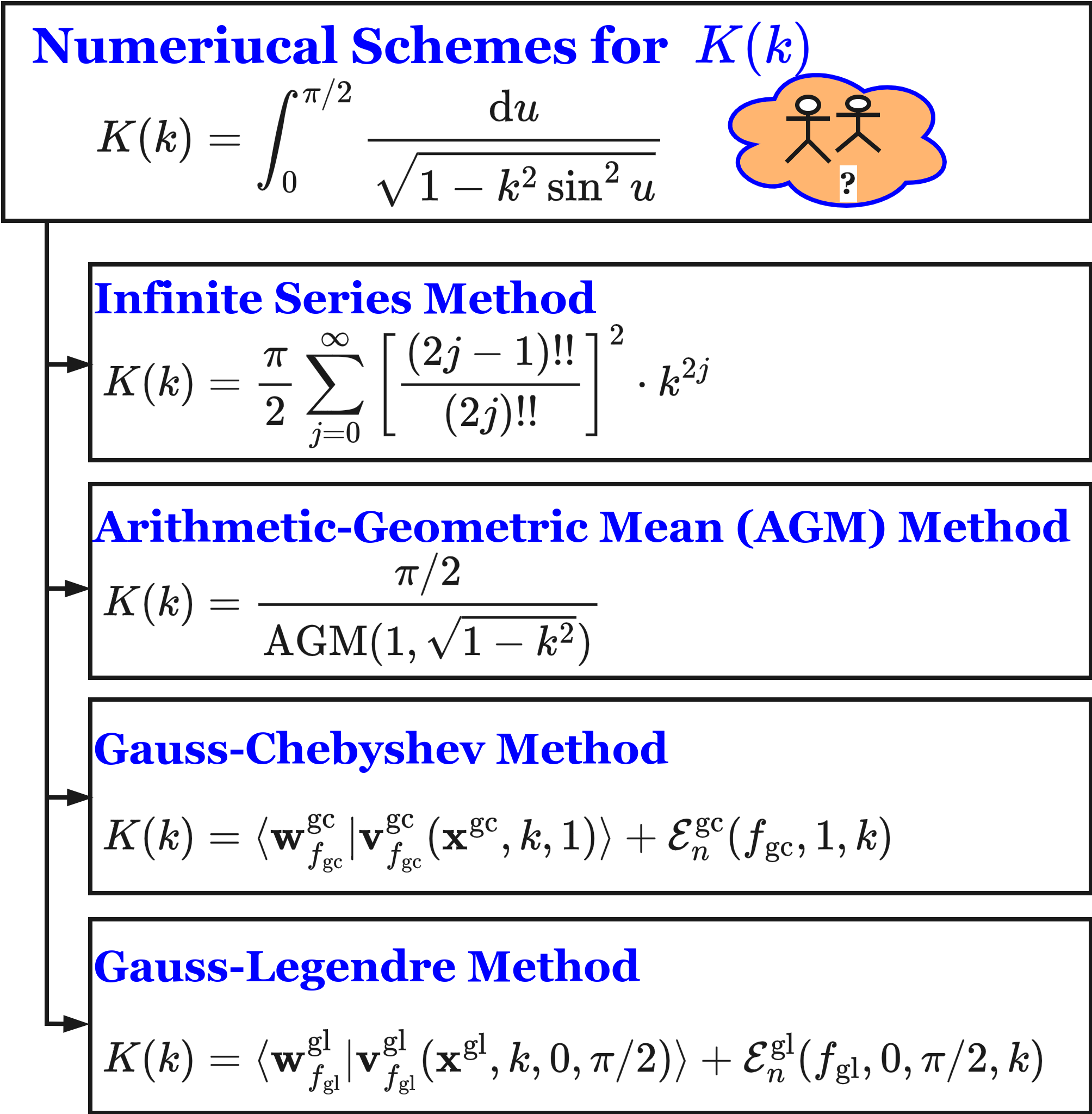} 
	\caption{Four methods for computing CEI-1}
	\label{fig-ways2K}
\end{figure}

In this paper, we will discuss four methods to calculate the special function $K(k)$ as shown in \Fig  \ref{fig-ways2K}.

The contents of this paper are organized as following: section \ref{sec-preliminaries} deals with the mathematical preliminaries for the infinite series methods, AGM method, Gauss-Chebyshev method and Gauss-Legendre method for computing CEI-1; section \ref{sec-numer-alg} copes with the numerical algorithms for the four methods mentioned above; section \ref{sec-vvt} gives the numerical results for verification and validation of the numerical solvers proposed and finally section \ref{sec-conclusions}  gives the conclusions. 

\section{Mathematical Preliminaries} \label{sec-preliminaries}

\subsection{Infinite Series Method}

The first method for computing the CEI-1 is to take the series expansion of variable $k$. Actually, with the help of the generalized Newton's binomial theorem, we have \cite{TAOCP1}
\begin{equation} \label{eq-alpha-series}
	(1+z)^\alpha = \sum^\infty_{j=0} \binom{\alpha}{j}z^j, \quad \abs{z}\le 1, \forall \alpha \in \Real
\end{equation}  
where 
\begin{equation}
	\binom{\alpha}{0} = 1
\end{equation}
and 
\begin{equation}
	\binom{\alpha}{j} = \frac{\alpha(\alpha-1)\cdots(\alpha-j+1)}{j!}=\prod^{j-1}_{t=0}\frac{\alpha-t}{j-t}
\end{equation}
is the binomial coefficient.
Put $z = -k^2 \sin^2 u$ and $\alpha = -\frac{1}{2}$, we can obtain
\begin{equation}
	\frac{1}{ \sqrt{1-k^2\sin^2u}} = \sum^\infty_{j=0}(-1)^j\binom{-\frac{1}{2}}{j}k^{2j}\sin^{2j}u.
\end{equation}
According to the identity
\begin{equation}
	(-1)^j \binom{-r}{j} = \binom{r+j-1}{j}, \quad r\in \Real, j\in \mathbb{N}
\end{equation}
and integral formula 
\begin{equation}
	\int^{\pi/2}_0 \sin^{2j} u \dif u = \frac{(2j-1)!!}{(2j)!!}\cdot\frac{\pi}{2}, \quad j\in \mathbb{N}
\end{equation}
where
\begin{equation}
\left\{
	\begin{split}
		&(2j-1)!! = (2j-1)(2j-3)\cdots 3 \cdot 1\\
		&(2j)!! = (2j)(2j-2)\cdots 4\cdot 2\\
		&0!! = (-1)!! = 1
	\end{split}
\right.
\end{equation}
are the fundamental results for the double factorial. Therefore, we can find the infinite series for the CEI-1, viz.,
\begin{equation} \label{eq-K-infinite-series}
K(k) = \frac{\pi}{2} \cdot  \sum^\infty_{j=0}\left[\frac{(2j-1)!!}{(2j)!!} \right]^2k^{2j} 
\end{equation}
Simple algebraic operations show that
\begin{equation} \label{eq-K-expansion}
K(k)=  \frac{\pi}{2} \cdot\left[ 1 + \frac{1}{4} k^2 + \frac{9}{64} k^4 + \frac{25}{256} k^6 + \cdots \right]
\end{equation}
Let
\begin{equation} \label{eq-cj}
	\quad c_j = \left[\frac{(2j-1)!!}{(2j)!!}\right]^2, \quad j\in \mathbb{N},
\end{equation}
then it is easy to find that
\begin{equation} \label{eq-K-series-k}
\left\{
\begin{aligned}
K(k) &= \frac{\pi}{2} \sum^\infty_{j=0}c_j k^{2j} \\
\scrd{\mathcal{K}}{cs}(m) &= \frac{\pi}{2} \sum^\infty_{j=0}c_j m^j
\end{aligned}
\right.
\end{equation}
according to \eqref{eq-km} and \eqref{eq-K-series-k}.

From a computational point of view, we just need to calculate the following convergent infinite series
\begin{equation} \label{eq-series}
	S(x) = \sum^\infty_{j=0} c_j x^j, \quad  x \in [0,1].
\end{equation} 
Therefore, the key for determining CEI-1 lies in designing an algorithm for computing the infinite series \eqref{eq-series} for $x=m$ or $x=k^2$ and designing a procedure to calculate the coefficient $c_j$ specified by \eqref{eq-cj}.

\subsection{AGM Method}

The second method for calculating the CEI-1 is the AGM method proposed by King \cite{King1924AGM}. The AGM method was first discovered by Lagrange and then rediscovered and pushed forward greatly by Gauss. The relation of AGM method and CEI-1 was also discussed by Borwein \cite{Borwein1987} for the purpose of calculating the constant $\pi$.   For the 2-dim sequence $\set{(u_i, v_i)}^\infty_{j=0}$
such that 
\begin{equation} \label{eq-AGM}
	\left\{
	\begin{split}
		u_{i+1} &= \frac{u_i +v_i}{2},\quad i\ge 0\\
		v_{i+1} & = \sqrt{u_i \cdot v_i}, \quad i\ge 0\\
		u_0 & = a \ge 0\\
		v_0 & = b \ge 0
	\end{split}
	\right.
\end{equation} 
it is easy to show that: both the limit of arithmetic mean $u_i$ and the limit of geometric mean $v_i$ exist, and the limits for $u_i$ and $v_i$ are equal. Obviously, the common limit just depends on the initial non-negative values  $a$ and $b$. This limit specified by the arithmetic-geometric mean (AGM) can be denoted by 
\begin{equation} \label{eq-agm-limit}
	\AGM(a,b) = \lim_{i\to \infty} u_i = \lim_{i\to \infty} v_i.
\end{equation}

The updating relation in equation \eqref{eq-AGM} implies that both the symmetric property
\begin{equation} \label{eq-AGM-symmetric}
	\AGM(a,b) = \AGM(b,a)
\end{equation}
and the iterative property
\begin{equation} \label{eq-AGM-iterative}
	\AGM(a,b) = \AGM\left(\frac{a+b}{2}, \sqrt{ab}\right)
\end{equation}
hold obviously.

There are lots of formula about the AGM and the CEI-1, we just list the following necessary ones according to our requirements:
\begin{align}
	\frac{1}{\AGM(1+x,1-x)} &= \sum^\infty_{j=0} \left[\frac{(2j-1)!!}{(2j)!!} \right]^2x^{2j} \label{eq-AGM-series}\\
	\frac{1}{\AGM(a,b)}
	&=\frac{2}{\pi}\int^{\frac{\pi}{2}}_0 \frac{\dif \theta}{\sqrt{a^2\cos^2 \theta + b^2\sin^2\theta}} \label{eq-AGM-integral}\\
	&=\frac{4}{(a+b)\pi}K\left(\abs{\frac{a-b}{a+b}}\right) \label{eq-AGM-K}
\end{align}
For $k\in [0,1], a = 1+k$ and $b = 1-k$, we can deduce that
\begin{equation}
	K(k) = \frac{\pi/2}{\AGM(1+k,1-k)}
\end{equation}
from \eqref{eq-AGM-series}. Consequently, we can obtain 
\begin{equation}
	\AGM(1+k,1-k) = \AGM(1, \sqrt{1-k^2})
\end{equation}
by the iterative property \eqref{eq-AGM-iterative}. Hence we have
\begin{equation}
K(k) =  \frac{\pi/2}{\AGM(1,\sqrt{1-k^2})}
\end{equation}

This is the famous AGM method proposed by King in 1924. 
Thus, in the sense of computation, the calculation of CEI-1 is equivalent to computing the function $\AGM(a,b)$ via iterative algorithm and calling it with 2-dim parameter $(a, b) = (1, \sqrt{1-k^2})$ or equivalently $(a,b) = (1+k, 1-k)$. Moreover, the symmetric property of AGM says that the position of $a$ and $b$ in $\AGM(a,b)$ can be swapped without changing the value.

\subsection{Gauss-Chebyshev Method} 

According to the Appendix \ref{appendix-ErrEstGC}, for the abstract function 
\begin{equation}
	\scrd{F}{C}(k)  = \int^a_{-a} \frac{f(x, k)}{\sqrt{a^2 - x^2}}\dif x 
	= \int^1_{-1} \frac{f(at, k)}{\sqrt{1-t^2}}\dif t,
\end{equation}
of the Chebyshev type  defined by the integral of $\displaystyle{\frac{f(x,k)}{\sqrt{a^2-x^2}}}$ with parameter $k$, we can computed 
it by the extended Gauss-Chebyshev formula of numerical integration. Let 
$\scru{x}{gc}_i$ denote the $i$-th root of $T_n(x) = \cos(n \arccos x)$, i.e. the $n$-th order Chebyshev polynomial of the first kind, and $\scru{w}{gc}_i$ denote the $i$-th weight corresponding to $\scru{x}{gc}_i$,  then we have
\begin{equation} \label{eq-x-w-gc}  
	\scru{x}{gc}_i = \cos \frac{(2i-1)\pi}{2n}, \quad \scru{w}{gc}_i = \frac{\pi}{n}, \quad 1\le i\le n
\end{equation}
and 
\begin{equation} \label{eq-F-Chebyshev-sum}
	\scrd{F}{C}(k) = \sum^{n}_{i=1} \scru{w}{gc}_i f(a\scru{x}{gc}_i, k) + \scru{\mathcal{E}}{gc}_n(f,a,k).
\end{equation}
where 
\begin{equation}
	\scru{\mathcal{E}}{gl}_n(f,a,k) = \frac{2\pi\cdot a^{2n}}{2^{2n}(2n)!}f^{(2n)}(\xi,k), \quad \xi \in (-a,a)
\end{equation}
is the approximation error.

We now introduce three $n$-dim vectors
\begin{equation}
	\left\{
	\begin{split}
		&\scru{\vec{x}}{gc} = \trsp{[\scru{x}{gc}_1, \cdots, \scru{x}{gc}_n]},\\
		&\scru{\vec{w}}{gc} = \trsp{[\scru{w}{gc}_1, \cdots, \scru{w}{gc}_n]},\\
		&\scru{\vec{v}}{gc}_{f}(k,a) = \trsp{[f(a\scru{x}{gc}_1,k), \cdots, f(a\scru{x}{gc}_n,k)]}
	\end{split}
	\right.
\end{equation}
then CEI-1 can be represented by the inner product of $\scru{\vec{w}}{gc}$ and $\scru{\vec{v}}{gc}_f(k)$, viz.,

\begin{equation} \label{eq-F-Chebyshev-innerproduct}
	\scrd{F}{C}(k) = \braket{\scru{\vec{w}}{gc}}{\scru{\vec{v}}{gc}_{f}(k,a)} + \scru{\mathcal{E}}{gc}_n(f,a,k)
\end{equation}
in which the inner product is defined by
\begin{equation}
	\braket{\vec{w}}{\vec{v}} = \sum^n_{i=1}w_iv_i, \quad \forall \vec{w}, \vec{v}\in \ES{R}{n}{1}.
\end{equation}

Our objective here is to give a computational method for CEI-1 by establishing the relation of CEI-1 and $\scrd{F}{C}(k)$. Let $t= \sin u$, then we have
\begin{equation}
	\begin{split}
		K(k) &= \int^{\pi/2}_0 \frac{\dif t}{\sqrt{1-k^2\sin^2 u}} 
		=\int^1_0 \frac{\dif t}{\sqrt{(1-t^2)(1-k^2t^2)}}\\
		&=\int^1_{-1}\frac{1}{\sqrt{1-t^2}}\cdot \frac{1/2}{\sqrt{1-k^2t^2}}\dif t
	\end{split}
\end{equation}
Now put $a=1$ and 
\begin{equation}\label{eq-f-gc}
	\scrd{f}{gc}(t,k) = \frac{1/2}{\sqrt{1-k^2t^2}}, \quad t \in(-1,1), k\in[0,\sqrt{2}/2) 
\end{equation}
then
\begin{equation}\label{eq-K-gc}
		\begin{split}
			K(k) &= \sum^{n}_{i=1} \scru{w}{gc}_i \scrd{f}{gc}(\scru{x}{gc}_i, k) + \scru{R}{gc}_n(\scrd{f}{gc})\\
			&=\braket{\scru{\vec{w}}{gc}}{\scru{\vec{v}}{gc}_{\scrd{f}{gc}}(\scru{\vec{x}}{gc},k,1)} + 
			\scru{\mathcal{E}}{gc}_n(\scrd{f}{gc},1,k)
		\end{split}
\end{equation}
Therefore, for any $\varepsilon >0$, the CEI-1 can be computed directly with the inner product of vectors $\scru{\vec{w}}{gc}$ and $\scru{\vec{v}}{gc}_{\scrd{f}{gc}}(k,1)$ if we take a sufficiently large $n$ such that $\abs{\scru{\mathcal{E}}{gc}_n(\scrd{f}{gc},1,k)}<\varepsilon$ for approximating the integral with fixed parameter $k$.

\subsection{Gauss-Legendre Method}

According to the Appendix \ref{appendix-ErrEstGL}, the abstract function  
\begin{equation}
		\scrd{F}{L}(k)
		=\int^b_a f(x,k)\dif x 
		= \frac{b-a}{2}\int^1_{-1} f\left(\frac{b+a}{2}+\frac{b-a}{2}t \right)\dif t
\end{equation}
defined by the integral of $f(x,k)$ with parameter $k$ such that $f(\cdot,k)\in C^{2n}[a,b]$ can be computed by the extended Gauss-Legendre formula of numerical integration. Let $\scru{x}{gl}_i$ denote the $i$-th root of the $n$-th order Legendre polynomial $P_n(x)$ and 
$\scru{w}{gl}_i$ denote the $i$-th weight corresponding to $\scru{x}{gl}_i$ such that
\begin{equation} \label{eq-weight-Legendre}
		\scru{w}{gl}_i 
		= \frac{2}{n}\cdot \frac{1}{P_{n-1}(\scru{x}{gl}_i)P'_n(\scru{x}{gl}_i)} 
		= \frac{2}{[1-(\scru{x}{gl}_i)^2]\cdot[P'_n(\scru{x}{gl}_i)]^2}
\end{equation}
where $P'_n(x)$ is the derivative function of $P_n(x)$, then we have
\begin{equation}
\begin{split}
		\scrd{F}{L}(k) 
		 &= \frac{b-a}{2}\sum^n_{i=1} \scru{w}{gl}_i 
		f\left(\frac{b+a+ (b-a)\scru{x}{gl}_i}{2}, k\right) \\ 
		 &\quad  + \scru{\mathcal{E}}{gl}_n(f,a,b,k)
\end{split}
\end{equation}
where 
\begin{equation}
	\scru{\mathcal{E}}{gl}_n(f,a,b,k) 
	= \frac{(n!)^4}{[(2n)!]^3}\cdot \frac{(b-a)^{2n+1}}{2n+1}f^{(2n)}(\xi,k)
\end{equation}
is the approximation error in which $\xi\in(a,b)$.

We remark that there are no analytical solutions for the roots of $P_n(x) = 0$ if $n \ge 4$. Fortunately, it can be solved by the Newton's method for nonlinear equation. For small $n$, say $n \le 20$, we can look up some available handbooks of mathematics to get the roots $\set{\scru{x}{gl}_i: 1\le i \le n}$. However, this method is not suggested if the $n$ is unknown or not sure. In this case, a feasible, convenient and effective way is to calculate the roots of $P_n(x)$ automatically with proper algorithms implemented via some concrete computer programming language and programs.

Let 
\begin{equation}
	\left\{
	\begin{split}
		&\scru{\vec{x}}{gl}  = \trsp{\left[\scru{x}{gl}_1, \cdots, \scru{x}{gl}_n\right]}\\
		&\scru{\vec{w}}{gl}  = \trsp{\left[\scru{w}{gl}_1, \cdots, \scru{w}{gl}_n\right]}\\
		&\scru{u}{gl}_i  = \frac{b+a + (b-a)\scru{x}{gl}_i}{2}, \quad 1\le i\le n \\
		&\scru{\vec{v}}{gl}_f(\scru{\vec{x}}{gl},k,a,b)  = \frac{b-a}{2}\trsp{\left[f(\scru{u}{gl}_i,k), \cdots, f(\scru{u}{gl}_n)\right]}
	\end{split}
	\right.
\end{equation}
then $\scrd{F}{L}(k)$ can be represented by the inner product of $\scru{\vec{w}}{gl}$ and $\scru{\vec{v}}{gl}_f(\scru{\vec{x}}{gl},k,a,b)$, viz.,
\begin{equation}
	\scrd{F}{L}(k) = \braket{\scru{\vec{w}}{gl}}{\scru{\vec{v}}{gl}_f(\scru{\vec{x}}{gl},k,a,b)} + \scru{\mathcal{E}}{gl}_n(f,a,b,k)
\end{equation}

Suppose that the root-weight pairs $(\scru{x}{gl}_i, \scru{w}{gl}_i)$ for $1\le i\le n$  were known, we can establish the relation of $\scrd{F}{L}(k)$ and CEI-1. For $k\in [0,1)$, let $a=0$, $b= \pi/2$, and 
\begin{equation} \label{eq-f-gl}
	\scrd{f}{gl}(u,k) = \frac{1}{\sqrt{1-k^2\sin^2u}}
\end{equation}
then 
\begin{equation} \label{eq-K-gl}
	\begin{split}
			K(k) &= \frac{\pi}{4}\sum^n_{i=1} \scru{w}{gl}_i \scrd{f}{gl}\left(\frac{\pi}{4}(1+\scru{x}{gl}_i), k\right) 
			+  \scru{\mathcal{E}}{gl}_n\left(\scrd{f}{gl},0,\frac{\pi}{2},k\right) \\
			&= \braket{\scru{\vec{w}}{gl}}{\scru{\vec{v}}{gl}_{\scrd{f}{gl}}(\scru{\vec{x}}{gl},k,0,\frac{\pi}{2})} +   \scru{\mathcal{E}}{gl}_n\left(\scrd{f}{gl},0,\frac{\pi}{2},k\right)
		\end{split}
\end{equation}
For any positive $\varepsilon >0$ and fixed parameter $k\in[0,1)$, we can choose a sufficiently large integer $n$ such that the approximation error is small enough, i.e., 
$$\abs{\scru{\mathcal{E}}{gl}_n\left(\scrd{f}{gl},0,\frac{\pi}{2},k\right)}
< \varepsilon.$$ 
The expression above for CEI-1 means that the essential work left for calculating $K(k)$ is to find the root-weight pairs for the Legendre polynomial $P_n(x)$. 

Although there are several equivalent definitions for the Legendre polynomial, the best way for computing the Legendre polynomial $P_n(x)$ is the following iterative formula
\begin{equation}
	\left\{
	\begin{split}
		P_n(x) & = \frac{2n-1}{n}xP_{n-1}(x) -\frac{n-1}{n}P_{n-2}(x), ~n \ge 2 \\
		P_0(x) & = 1; \\
		P_1(x) & = x.
	\end{split}
	\right.
\end{equation}
This formula implies that we can compute $P_n(x)$  one by one with an iterative process.

The Newton's method for solving root of a nonlinear equation $f(x,\alpha) = 0$ with a constant parameter $\alpha$ is specified by the iterative process with the form
\begin{equation} \label{eq-Newton-root-iter}
	\begin{split}
		x[j+1] &= \scrd{\mathcal{A}}{Newton}(x[j])\\
		&= x[j] - \frac{f(x[j],\alpha)}{f'(x[j],\alpha)}, \quad j = 0, 1, 2, \cdots
	\end{split}
\end{equation}
which starts from the initial value (guess of the root) $x[0] = x_0$ and ends when 
\begin{equation}
	\abs{x[j+1]-x[j]} < \varepsilon
\end{equation}
for sufficiently large $j$ and fixed computational error $\varepsilon$. The updating mapping $\scrd{\mathcal{A}}{Newton}: x[j]\mapsto x[j+1]$ is a contractive mapping defined by
\begin{equation} \label{eq-Newton-map}
	\begin{split}
		\scrd{\mathcal{A}}{Newton}: \Real &\to \Real \\
		x &\mapsto x - \frac{f(x,\alpha)}{f'(x,\alpha)} 
	\end{split}
\end{equation}

For the Legendre polynomial $P_n(x)$, there are $n$ different roots $\scru{x}{gl}_i$ where $i\in \set{1, 2, \cdots, n}$. The iterative process can be described by
\begin{equation} \label{eq-init-values-for-Pn}
	\left\{
	\begin{split}
		\scru{x}{gl}_i[j+1] &= \scru{x}{gl}_i[j] - \frac{P_n(\scru{x}{gl}_i[j])}{P'_n(\scru{x}{gl}_i[j])}, \quad j = 0, 1, 2 \cdots \\
		\scru{x}{gl}_i[0]   &= \cos \frac{\pi(i + 0.75)}{n + 0.25}, \quad 1\le i \le n. 
	\end{split}
	\right.
\end{equation}
This iterative process is the basis for designing algorithm for solving the roots of interest. Note that we can set $f(x,\alpha) = P_n(x)$ where $\alpha = n$ here. 

\section{Numerical Algorithms} \label{sec-numer-alg}

Different numerical algorithms can be obtained according to the four methods discussed to compute the CEI-1 efficiently, robustly and precisely. In order to broaden the applications of the algorithms, we will generalize some specific constants to corresponding variables.

\subsection{Infinite Series Algorithm}

The equation \eqref{eq-K-infinite-series} and its connection with the general infinite series show that there are some key  issues for computing CEI-1:
\begin{itemize}
	\item designing an algorithm to compute the general infinite series $S(x) = \sum\limits_{j=0}^\infty c_j x^j$ which consists of four sub-tasks:
	\begin{itemize}
		\item calculating the coefficient $c_j$ for each $j\in \mathbb{N}$,
		\item computing the value of the general term $c_j x^j$ for fixed $x\in \Real$,
		\item deciding whether the general term $c_jx^j$ is negligible  or not for the sum of the infinite series, and 
		\item calculating the sum of the not-negligible terms;
	\end{itemize}
	\item specifying the parameters involved 
	\begin{itemize}
		\item designing a procedure to compute the coefficient $c_j$ in the infinite series of CEI-1, and
		\item specifying the variable $x$ according to its connection with the variable $k$, i.e., $x = k^2$.
	\end{itemize}
	\item calling the procedure for the infinite series and multiplying some constant if necessary.
\end{itemize}

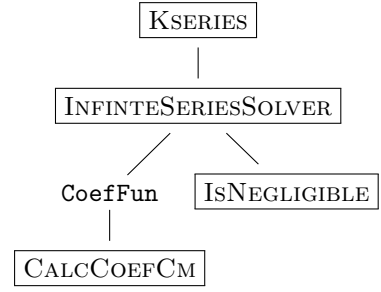
\begin{figure}[h]
	\centering
	\begin{forest}
		[\fbox{\ProcName{Kseries}}
		[\fbox{\ProcName{InfinteSeriesSolver}}
		[\PtrToFun{CoefFun}[\fbox{\ProcName{CalcCoefCm}}]]
		[\fbox{\ProcName{IsNegligible}}]
		]
		] 
	\end{forest}
	\caption{Procedures involved in calculating the complete elliptic integral CEI-1 with infinite series method.}
	\label{fig-Kseries-procedures}
\end{figure}

The procedures for calculating CEI-1 with infinite series method is shown in 
\Fig \ref{fig-Kseries-procedures}, in which \PtrToFun{CoefFun} is a procedure name (pointer to function or similar entity)  used to calling the procedure for computing $c_j$. 

There is a practical technique for computing $c_j$. The double factorials $(2j)!!$ and $(2j-1)!!$ may be large numbers if $j>10$ and the overflow problem may happen in some computer language if the type of data is integer. However, we can cope with the the value of the fraction safely instead of dealing with its denominator and numerator separately. Actually, we have

\begin{equation}
		\frac{(2j-1)!!}{(2j)!!}
		=\prod^{j-1}_{t=0}\frac{2(j-t)-1}{2(j-t)}
		=\prod^{j}_{t=1}\left[
		1-\frac{0.5}{t}\right], \quad t \ge 1.
\end{equation}

Thus the fraction $\frac{(2j-1)!!}{(2j)!!} $ can be computed by multiplying  $1-0.5/t$ one by one with a simple iteration, which not only reduces the computational complexity but also avoid the overflow problem of the factorial operations. \Algr \ref{alg-CalcCoefCm} gives the steps of procedure \ProcName{CalcCoefCm}  for computing the coefficient $c_m =\frac{(2m-1)!!}{(2m)!!}$  by the for-loop construction in computer programming language.

\begin{breakablealgorithm}
	\caption{Calculate the coefficient $c_j=\left(\frac{(2j-1)!!}{(2j)!!}\right)^2$ for the $j$-th term of the infinite series for CEI-1} \label{alg-CalcCoefCm}
	\begin{algorithmic}[1]
		\Require Non-negative integer $j$;
		\Ensure The coefficient $c_j$ of the infinite series in the function $K(k)=\dfrac{\pi}{2}\sum\limits^\infty_{j=0}c_jk^{2j}$
		\Function{CalcCoefCm}{$j$}
		\If{$j == 0$}
		\State \Return 1;
		\EndIf
		\State $\cpvar{prod} \gets 1$;
		\For{$t\in\seq{1,2,\cdots,j-1,j}$}
		\State $\cpvar{prod} \gets \cpvar{prod} \cdot (1 - 0.5/t)$;
		\EndFor
		\State $c_j \gets \cpvar{prod} \cdot \cpvar{prod}$;
		\State \Return $c_j$;
		\EndFunction
	\end{algorithmic}
\end{breakablealgorithm}

For a fixed $x$, the value of the infinite series $S(x) = \sum\limits^\infty_{j=0}c_mx^m$ is defined by
\begin{equation}
	S(x) = \lim_{i\to \infty} S_i(x), \quad S_i(x)=\sum^{i-1}_{j=0}c_jx^j
\end{equation}
The Cauchy's criteria of convergence shows that the sequence $\set{S_m(x)}^\infty_{j=0}$ converges if and only if 
\begin{equation}
	\forall \varepsilon>0, \quad \abs{S_i(x) - S_{i-1}(x)} = \abs{c_ix^i} < \varepsilon.
\end{equation}
In practice, the relative error relation $\frac{\abs{S_i(x) - S_{i-1}(x)} }{\abs{S_i(x)}}< \varepsilon$ is better than the absolute error relation $\abs{S_i(x) - S_{i-1}(x)}< \varepsilon $.  This fact can be used to design the procedure for deciding whether the general term $c_ix^i$ is negligible or not.

\begin{breakablealgorithm}
	\caption{Decide whether the general term $c_jx^j$ is negligible or not}  \label{alg-IsNotNeglegible}
	\begin{algorithmic}[1]
		\Require The value of the general term $c_jx^j$: $\cpvar{term}$, the sum of finite terms: $\cpvar{sum}$, and the value of error (precision) $\varepsilon$ with default value $\varepsilon = 10^{-9}$.
		\Ensure  The logical value $\True$ or $\False$
		\Function{IsNegligible}{$\cpvar{term},\cpvar{sum},\varepsilon$}
		\State Set a sufficiently small number: $\delta \gets 10^{-10}$;
		\If{$\abs{\cpvar{term}}/(\abs{\cpvar{sum}} + \delta) < \varepsilon$}
		\State   \Return $\True$;
		\Else
		\State  \Return $\False$; 
		\EndIf
		\EndFunction
	\end{algorithmic}
\end{breakablealgorithm}  

\Algr \ref{alg-IsNotNeglegible}, named by \ProcName{IsNegligible}, is utilized to decide whether the $j$-th term $c_j x^j$ is important for the sum of infinite series $\sum\limits_{j=0}^\infty c_jx^j$ or not. We remark that in \Algr \ref{alg-IsNotNeglegible} the small number $\delta$ is used to improve the robustness when $\abs{S_j(x)}$ is near to zero.

The procedure for solving the infinite series $\sum\limits^\infty_{j=0}c_jx^j$ can be designed with a do-while or while-do loop, which is named with \ProcName{InfiniteSeriesSolver}  and shown in \Algr \ref{alg-InfinteSeriesSolver}. Obviously, the computation of CEI-1 is treated as a special problem of sovling the infinite series by designing a system with the strategy of top-down. 

\begin{breakablealgorithm} 
	\caption{Calculate the sum of the convergent infinite series $\displaystyle\sum^\infty_{j=0} c_j x^j$ (\ProcName{InfinteSeriesSolver})}
	\label{alg-InfinteSeriesSolver}
	\begin{algorithmic}[1]
		\Require Pointer to function $\PtrToFun{CoefFun}: \mathbb{N}\to \Real$  for computing the coefficient $c_j$,  variable $x$, and precision $\varepsilon$ with default value $\varepsilon = 10^{-9}$.
		\Ensure The value of $\displaystyle \sum^\infty_{j=0} c_j x^j$
		\Function{InfinteSeriesSolver}{$\PtrToFun{CoefFun}, x, \varepsilon$}
		\State Set initial values: $\cpvar{sum}\gets 0, \cpvar{term} \gets 0,
		j \gets 0$;
		\Repeat
		\State $c_j \gets \PtrToFun{CoefFun}(j)$;
		\State $\cpvar{term} \gets c_j\cdot x^j$;
		\State $\cpvar{sum} \gets \cpvar{sum} + \cpvar{term}$;
		\State $j \gets j + 1$;
		\Until{$(~!~\ProcName{IsNeglibible}(\cpvar{term}, \cpvar{sum}, \varepsilon))$}
		\State \Return \cpvar{sum};
		\EndFunction
	\end{algorithmic}
\end{breakablealgorithm}

It should be noted that the function $\PtrToFun{CoefFun}: \mathbb{N}\to \Real$ is used as an argument for the procedure \ProcName{InfiniteSeries} in \Algr \ref{alg-InfinteSeriesSolver}, which implies that the \ProcName{InfiniteSeries} is a high order function and the pointer to function will be involved for the implementation with different computer programming languages. We remark that the pointer to function is a concept in the C programming language and its counterparts in Lisp/Python/Java/C++ and MATLAB/Octave are lambda expression and function handle respectively.

Now we can design the procedure  for calculating the CEI-1  conveniently. \Algr \ref{alg-Kseries} describes the procedure \ProcName{Kseries} for computing CEI-1 with the infinite series method. 
\begin{breakablealgorithm} 
	\caption{Calcualting the function CEI-1 with the infinite series method (\ProcName{Kseries})}
	\label{alg-Kseries}
	\begin{algorithmic}[1]
		\Require The variable $k\in[0,1)$, precision $\varepsilon$ with the default value $\varepsilon = 10^{-9}$.
		\Ensure The value of CEI-1 
		\Function{Kseries}{$k, \varepsilon$}
		\State $\PtrToFun{CoefFun}\gets \ProcName{CalcCoefCm}$;
		\State $x\gets k^2$; // Attention, please!
		\State $\cpvar{sum} \gets \ProcName{InfinteSeriesSolver}(\PtrToFun{CoefFun}, x, \varepsilon)$;
		\State $K \gets \cpvar{sum} \cdot \pi/2$;
		\State \Return $K$;
		\EndFunction
	\end{algorithmic}
\end{breakablealgorithm}

\subsection{AGM Algorithm}

\begin{figure}[h] 
	\centering
	\begin{forest}
		[\fbox{\ProcName{Kagm}}
		[\fbox{\ProcName{AGM}}
		[\fbox{\ProcName{FixedPointSolver}}
		[\PtrToFun{Update}[\fbox{\ProcName{AgmUpdate}}]]
		[\PtrToFun{Dist}[\fbox{\ProcName{DistEuclid}}]]
		]
		[$a$]
		[$b$]
		[$\varepsilon$]
		]
		] 
	\end{forest}
	\caption{Procedures for calculating CEI-1 with AGM method.}
	\label{fig-Kagm-procedures}
\end{figure}
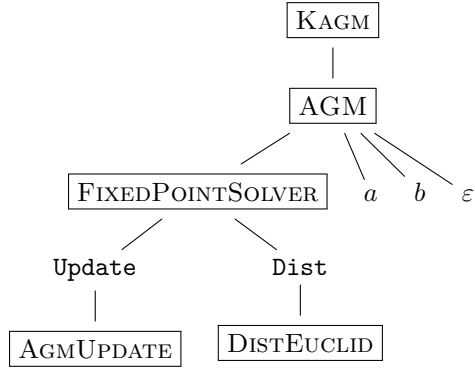

The key issue of the AGM is to compute the limits 
$\lim\limits_{i\to \infty} u_i$ and $\lim\limits_{i\to \infty} v_i$ according to \eqref{eq-AGM}. 
The procedures for calculating CEI-1 with the AGM method is shown in \Fig \ref{fig-Kagm-procedures}, in which \PtrToFun{Update} and \PtrToFun{Dist} are two pointers to functions and will be assigned by \ProcName{AgmUpdate} in \Algr \ref{alg-AgmUpdate} and \ProcName{DistEuclid} in \Algr \ref{alg-DistEuclid} respectively.

Let $\vec{x} = \trsp{[x_1, x_2]}\in \ES{R}{2}{1}$, 
then the AGM method can be converted to a 2-dim fixed-point problem with the iterative formula
\begin{equation} \label{eq-fixed-point}
	\vec{x}_{i+1} = \scrd{\mathcal{A}}{agm}(\vec{x}_i), \quad i= 0, 1, 2, \cdots
\end{equation}
such that $\vec{x}_i \in \ES{R}{2}{1}_+$ and 
\begin{equation} \label{eq-agm-iter}
	\scrd{\mathcal{A}}{agm}(\vec{x}) = \begin{bmatrix}
		(x_1+x_2)/2\\
		\sqrt{x_1x_2}
	\end{bmatrix}, \quad 
	\vec{x}_0 = \begin{bmatrix}a\\b\end{bmatrix}
\end{equation}
for non-negative $a$ and $b$. 
The Banach's fixed-point theorem in functional analysis shows that \eqref{eq-fixed-point} is convergent when $\scrd{\mathcal{A}}{agm}$ is contractive ($\norm{\scrd{\mathcal{A}}{agm}} < 1$). Without doubt, \eqref{eq-agm-iter} satisfies the convergent condition for any non-negative components $a$ and $b$ of the initial vector $\vec{x}_0$. 

The updating process in AGM method gives a typical nonlinear mapping which is described by the procedure \ProcName{AgmUpdate} in \Algr \ref{alg-AgmUpdate}.
\begin{breakablealgorithm}
	\caption{AGM updating (\ProcName{AgmUpdate})}
	\label{alg-AgmUpdate}
	\begin{algorithmic}[1]
		\Require {2-dim vector $\vec{x}=\trsp{[x_1, x_2]}$}
		\Ensure 2-dim vector $\vec{y} = \mathcal{A}(\vec{x})$ 
		\Function{AgmUpdate}{$\vec{x}$}
		\State $y_1 \gets (x_1+x_2)/2$;
		\State $y_2 \gets \sqrt{x_1 x_2}$;
		\State $\vec{y} \gets \trsp{[y_1, y_2]}$;
		\State \Return $\vec{y}$;
		\EndFunction
	\end{algorithmic}
\end{breakablealgorithm}

It is wise to design a general algorithm for solving the fixed-point in a general way since the fixed-point may be frequently encountered in different STEM problems. The computation of fixed-point is based on the concept of norm. In consequence, it is necessary to design the norm for the vector $\vec{x}$. Usually, the Euclidean norm $\norm{\vec{x}} = \sqrt{\sum\limits^n_{j=1}x^2_j}$ for $\vec{x}=\trsp{[x_1, \cdots, x_n]}\in \ES{R}{n}{1}$ will be a good choice and the induced 2-dim Euclidean distance will be
\begin{equation}
	\dist(\vec{x}, \vec{y}) =\normp{\vec{x}-\vec{y}}{2} 
	= \sqrt{\sum^n_{j=1}(x_j-y_j)^2}
\end{equation}
for $\vec{x} = \trsp{[x_1, \cdots, x_n]}$ and $\vec{y} = \trsp{[y_1, \cdots, y_n]}$. 
The procedure \ProcName{DistEuclid} described in \Algr \ref{alg-DistEuclid} is used to compute the $n$-dim Euclidean distance.  
\begin{breakablealgorithm}
	\caption{Calculate the Euclidean distance of two $n$-dim vectors (\ProcName{DistEuclid}).}
	\label{alg-DistEuclid}
	\begin{algorithmic}[1]
		\Require Two $n$-dim vectors $\vec{x},\vec{y}\in \ES{R}{n}{1}$.
		\Ensure The distance between $\vec{x}$ and $\vec{y}$.
		\Function{DistEuclid}{$\vec{x}, \vec{y}, n$}
		\State $\cpvar{sum} \gets 0$ ;
		\For{$i\in \seq{1, \cdots, n}$}
		\State $\cpvar{sum} \gets \cpvar{sum} + (x_i - y_i)^2$;
		\EndFor
		\State $\cpvar{dist} \gets \sqrt{\cpvar{sum}}$;
		\State \Return $\cpvar{dist}$;
		\EndFunction
	\end{algorithmic}
\end{breakablealgorithm}

The iterative process of AGM implies a special kind of iteration for the fixed-point, which can be computed in a generic way with the help of while-do or do-while loop. The procedure \ProcName{FixedPointSolver}, as shown in \Algr \ref{alg-AbstractFixedPointSolver}, computes the $n$-dim fixed-point iteratively. Note that the pointers to functions, viz. \PtrToFun{Update} and \PtrToFun{Dist} should be assigned with concrete procedures according to the data type of $\vec{x}$. 

\begin{breakablealgorithm}
	\centering
	\caption{Compute the abstract fixed-point iteratively 
		(\ProcName{AbstractFixedPointSolver})}\label{alg-AbstractFixedPointSolver}
	\begin{algorithmic}[1] 
		\Require Dimension $n$, initial data $\vec{x}_0$, precision $\varepsilon$, updating function $\PtrToFun{Update}$ for the contractive mapping $\mathcal{A}: \mathscr{V}\to \mathscr{V}$, and distance function $\PtrToFun{Dist}: \mathscr{V}\times \mathscr{V}\to \Real^+$.
		\Ensure The fixed-point $\vec{x}$ such that $\vec{x} = \mathcal{A}(\vec{x})$.
		\Function{AbstractFixedPointSolver}{$n$, \PtrToFun{Update}, \PtrToFun{Dist}, $\vec{x}_0$, $\varepsilon$}
		\State $\cpvar{guess} \gets \vec{x}_0$;
		\State $\cpvar{improve}\gets  \PtrToFun{Update}(\vec{x}_0)$;
		\While{$(\PtrToFun{Dist}(\cpvar{guess}, \cpvar{improve}, n) \ge \varepsilon)$}
		\State $\cpvar{guess} \gets \cpvar{improve}$;
		\State $\cpvar{improve}\gets  \PtrToFun{Update}(\cpvar{guess})$;
		\EndWhile
		\State \Return \cpvar{improve};
		\EndFunction
	\end{algorithmic}
\end{breakablealgorithm}

If the general Banach space $\mathscr{V}$ is just the $n$-dim Euclidean space $\ES{R}{n}{1}$ and the distance is just the Euclidean distance, then the computation of distance will be easy. Consequently, we can design a procedure \ProcName{FixedPointSolver} for computing the fixed-point in $\ES{R}{n}{1}$. For more details, please see 
\Algr \ref{alg-FixedPointSolver}. 

For the AGM problem concerned, we have the dimension $n = 2$ and the 2-dim vectors can be represented by a 1-dim array of size $2$ or \cpvar{struct} with two components when coding with some concrete computer programming language. Therefore, the AGM method can be implemented by taking the
procedure \ProcName{FixedPointSolver}  in \Algr \ref{alg-FixedPointSolver}  and assigning the pointers to function \PtrToFun{Update} with the procedures \ProcName{AgmUpdate} in 
\Algr \ref{alg-AgmUpdate}. Once the iteration ends, we can return the arithmetic mean or the geometric mean in the 2-dim vector according to \eqref{eq-agm-limit}. For more details, please see the procedure \ProcName{AGM} in \Algr \ref{alg-AGM}.

\begin{breakablealgorithm}
	\centering
	\caption{Compute the $n$-dim fixed-point $\vec{x}\in \ES{R}{n}{1}$ such that $\vec{x} = \mathcal{A}(\vec{x})$ with iterative method  (\ProcName{FixedPointSolver}). }
	\label{alg-FixedPointSolver}
	\begin{algorithmic}[1]
		\Require Dimension $n$, initial data $\vec{x}_0$, precision $\varepsilon$, updating function $\PtrToFun{Update}$ for the contractive mapping $\mathcal{A}: \ES{R}{n}{1}\to \ES{R}{n}{1}$.
		\Ensure The $n$-dim fixed-point $\vec{x}\in \ES{R}{n}{1}$ such that $\vec{x} = \mathcal{A}(\vec{x})$.
		\Function{FixedPointSolver}{$n$, \PtrToFun{Update}, $\vec{x}_0$, $\varepsilon$}
		\State $\cpvar{guess} \gets \vec{x}_0$;
		\State $\cpvar{improve}\gets  \PtrToFun{Update}(\vec{x}_0)$;
		\While{$(\ProcName{DistEuclid}(\cpvar{guess}, \cpvar{improve}, n) \ge \varepsilon)$}
		\State $\cpvar{guess} \gets \cpvar{improve}$;
		\State $\cpvar{improve}\gets  \PtrToFun{Update}(\cpvar{guess})$;
		\EndWhile
		\State \Return \cpvar{improve};
		\EndFunction
	\end{algorithmic}
\end{breakablealgorithm}

\begin{breakablealgorithm}
	\caption{Calculate the value of the arithmetic-geometric mean of non-negative values $a$ and $b$ (\ProcName{AGM}).}
	\label{alg-AGM}
	\begin{algorithmic}[1]
		\Require Non-negative values $a$ and $b$ and precision $\varepsilon$.
		\Ensure The value of the AGM for $a$ and $b$.
		\Function{AGM}{$a,b, \varepsilon$}
		\State $\vec{x}_0 \gets \trsp{[a,b]}$; 
		\State $\PtrToFun{Update}\gets \ProcName{AgmUpdate}$;
		\State $\PtrToFun{Dist}\gets \ProcName{DistEuclid}$;
		\State $\cpvar{dim} \gets 2$;
		\State $\vec{x}\gets \ProcName{FixedPointSolver}
		(\cpvar{dim}, \PtrToFun{Update}, \PtrToFun{Dist}, \vec{x}_0, \varepsilon)$; 
		\State \Return $x_1$;  // return one component of $\vec{x}$;
		\EndFunction
	\end{algorithmic}
\end{breakablealgorithm}

The procedure for computing CEI-1 with the AGM method, named by \ProcName{Kagm}, is illustrated in \Algr~\ref{alg-Kagm} in which the procedure \ProcName{AGM} is called directly by setting proper parameters $a \gets 1$ and $b \gets \sqrt{1-k^2}$.  

\begin{breakablealgorithm}
	\caption{Calculate the CEI-1 with the AGM method (\ProcName{Kagm})}
	\label{alg-Kagm}
	\begin{algorithmic}[1]
		\Require The variable $k\in[0,1)$, precision $\varepsilon$ with default value $\varepsilon = 10^{-9}$.
		\Ensure The value of CEI-1
		\Function{Kagm}{$k, \varepsilon$}
		\State $a \gets 1$;
		\State $b \gets \sqrt{1-k^2}$;
		\State $K \gets (\pi/2)/\ProcName{AGM}(a,b, \varepsilon)$;
		\State \Return $K$;
		\EndFunction
	\end{algorithmic}
\end{breakablealgorithm}

\subsection{Gauss-Chebyshev Algorithm}

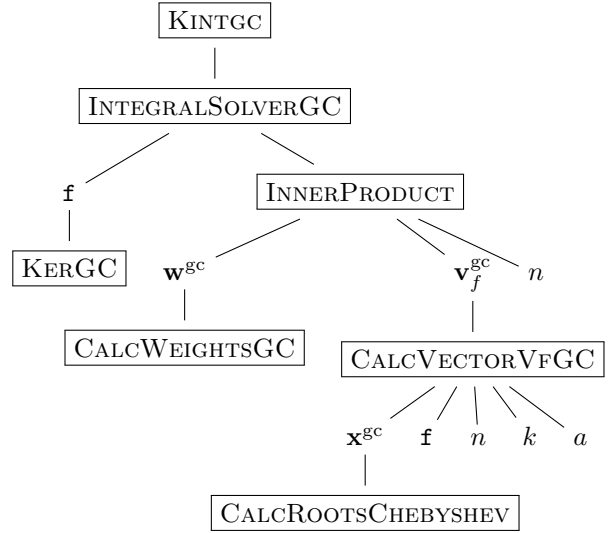
\begin{figure}[h] \label{fig-Kintgc-relations}
	\centering
	\begin{forest}
		[\fbox{\ProcName{Kintgc}}
		[\fbox{\ProcName{IntegralSolverGC}}
		[\PtrToFun{f}[\fbox{\ProcName{KerGC}}]]
		[\fbox{\ProcName{InnerProduct}}
		[$\scru{\vec{w}}{gc}$
		[\fbox{\ProcName{CalcWeightsGC}}]    	
		]
		[$\scru{\vec{v}}{gc}_f$
		[\fbox{\ProcName{CalcVectorVfGC}}
		[$\scru{\vec{x}}{gc}$
		[\fbox{\ProcName{CalcRootsChebyshev}}]
		]
		[\PtrToFun{f}]	
		[$n$]
		[$k$]	
		[$a$]
		]
		]
		[$n$]
		] 
		]	
		] 
	\end{forest}
	\caption{Procedures involved in calculating complete elliptic integral CEI-1 with Gauss-Chebyshev method.}
	\label{fig-Kintgc-procedures}
\end{figure}

The algorithm of Gauss-Chebyshev for computing CEI-1 is relatively simple due to the
equations \eqref{eq-x-w-gc}, \eqref{eq-f-gc} and \eqref{eq-K-gc}.
\Fig\ref{fig-Kintgc-procedures} demonstrates the procedures involved in calculating CEI-1 with the Gauss-Chebyshev method. 

The computation of CEI-1 is recognized as a special problem of solving the general numerical integration 
$\scrd{F}{C}(k)$ according to  \eqref{eq-F-Chebyshev-innerproduct}, which corresponds to the procedure 
\ProcName{IntegralSolverGC} in \Algr \ref{alg-IntegralSolverGC}, 
by setting the pointer to function \PtrToFun{f} with \ProcName{KerGC} in \Algr \ref{alg-KerGC}.

We now give some interpretation for designing  the algorithm: 
\begin{itemize}
	\item The procedure \ProcName{KerGC} in \Algr \ref{alg-KerGC} is defined to compute the function $\displaystyle \scrd{f}{gc}(t,k) = \frac{1/2}{\sqrt{1-k^2t^2}}$ and
	use a pointer to function \cpvar{KerFun} to represent the $\scrd{f}{gc}$ in the design and implementation of the algorithm. 
	\item 
	The root vector $\scru{\vec{x}}{gc}$ is set by the procedure \ProcName{CalcRootschebyshev} in  \Algr \ref{alg-CalcRootsChebyshev}. 
	\item 
	The value of integral is computed by the procedure \ProcName{InnerProduct} in \Algr \ref{alg-InnerProduct}.
	\item The arguments of \Algr \ref{alg-CalcRootsChebyshev}, namely the vector $\scru{\vec{w}}{gc}$ and vector $\scru{\vec{v}}{gc}_f(\scru{\vec{x}}{gc}, k, a)$,  are calculated by the procedures \ProcName{CalcWeightsGC} in \Algr \ref{alg-CalcWeightsGC} and 
	\ProcName{CalcVectorVfGC} in \Algr \ref{alg-CalcVectorVfGC} respectively. 
\end{itemize}

\begin{breakablealgorithm}
	\caption{Calcualte the integral  (abstract function)   
		\[
		\scrd{F}{C}(y) = \int^a_{-a} \frac{f(x,k)}{\sqrt{a^2 - x^2}}\dif x = \braket{\scru{\vec{w}}{gc}}{\scru{\vec{v}}{gc}_f(\scru{\vec{x}}{gc}, k,a)} + \mathcal{E}_n(f)
		\] 
		with the parameter $k$ via the $n$-nodes Gauss-Chebyshev numerical integration method (\ProcName{IntegralSolverGC}).}
	\label{alg-IntegralSolverGC}
	\begin{algorithmic}[1]
		\Require Variable $k$, pointer to function $\PtrToFun{f}: \Real \times \Real \to \Real$, positive integral $n$, $n$-dim root vector $\scru{\vec{x}}{gc}$, $n$-dim weight vector $\scru{\vec{w}}{gc}$, positive  $a$.
		\Ensure The value  of $\scrd{F}{C}(y)$.
		\Function{IntegralSolverGC}{$k, \PtrToFun{f}, n, \scru{\vec{x}}{gc}, \scru{\vec{w}}{gc}, a$}
		\State $\scru{\vec{v}}{gc}\gets \ProcName{CalcVectorVfGC}( 
		\PtrToFun{f}, n, \scru{\vec{x}}{gc}, k, a)$;
		\State $\cpvar{j}\gets \ProcName{InnerProduct}(\scru{\vec{w}}{gc}, \scru{\vec{v}}{gc}, n)$;
		\State \Return \cpvar{j};
		\EndFunction
	\end{algorithmic}
\end{breakablealgorithm}

\begin{breakablealgorithm}
	\caption{Compute the function $\scrd{f}{gc}(t,k)$ in the Gauss-Chebyshev method for calculating the value CEI-1 (\ProcName{KerGC})}
	\label{alg-KerGC}
	\begin{algorithmic}[1]
		\Require {Variables $t\in [-1,1], k\in [0,1)$}
		\Ensure The value of $\scrd{f}{gc}(t,k)$
		\Function{KerGC}{$t,k$}
		\State $u \gets kt$;
		\State $v\gets 0.5/\sqrt{1-u^2}$;
		\State \Return $v$;
		\EndFunction
	\end{algorithmic}
\end{breakablealgorithm}

\begin{breakablealgorithm}
	\caption{Calculate the root vector for the $n$ roots of the Chebyshev polynomial $T_n(x)=\cos(n\arccos x)$ (\ProcName{CalcRootsChebyshev})}
	\label{alg-CalcRootsChebyshev}
	\begin{algorithmic}[1]
		\Require {Positive integer $n$.}
		\Ensure The $n$ roots of $T_n(x)$.
		\Function{CalcRootsChebyshev}{$n$}
		\For{$i\in \set{1, \cdots, n}$}
		\State $\scru{x}{gc}_i \gets \cos \cfrac{(2i-1)\pi}{2n}$;
		\EndFor
		\State \Return $\scru{\vec{x}}{gc}$;
		\EndFunction
	\end{algorithmic}
\end{breakablealgorithm}

\begin{breakablealgorithm} 
	\caption{Compute the inner product of two $n$-dim vectors. (\ProcName{InnerProduct})}
	\label{alg-InnerProduct}
	\begin{algorithmic}[1]
		\Require Two vectors $\vec{w}$ and $\vec{v}$, dimension $n$ of the vectors.
		\Ensure Inner product $\braket{\vec{w}}{\vec{v}}=\sum\limits^n_{i=1}w_iv_i$.
		\Function{InnerProduct}{$\vec{w},\vec{v}, n$}
		\State $\cpvar{sum} \gets 0$;
		\For{$i\in \set{1,\cdots, n}$}
		\State $\cpvar{sum} \gets \cpvar{sum} + w_i\cdot v_i$;
		\EndFor
		\State \Return \cpvar{sum};
		\EndFunction
	\end{algorithmic}
\end{breakablealgorithm}

\begin{breakablealgorithm}
	\caption{Calculate the weight vector $\scru{\vec{w}}{gc}$ for the abstract function (parametric integral with parameter $k$)
		$\scrd{F}{C}(k) = \displaystyle{\int^a_{-a} \frac{f(x,k)}{\sqrt{a^2 -x^2}}\dif x} = \braket{\scru{\vec{w}}{gc}}{\scru{\vec{v}}{gc}_f(\scru{\vec{x}}{gc},k,a)} + \mathcal{E}_n(f)$ with the $n$-nodes Gauss-Chebyshev numerical integration method (\ProcName{CalcWeightsGC})}
	\label{alg-CalcWeightsGC}
	\begin{algorithmic}[1]
		\Require Positive integer $n$
		\Ensure The weight vector $\scru{\vec{w}}{gc}$
		\Function{CalcWeightsGC}{$n$}
		\For{$i\in \set{1, \cdots, n}$}
		\State $\scru{w}{gc}_i \gets \pi/n$;
		\EndFor
		\State \Return $\scru{\vec{w}}{gc}$;
		\EndFunction
	\end{algorithmic}
\end{breakablealgorithm}

\begin{breakablealgorithm}
	\caption{Calculate the vector $\scru{\vec{v}}{gc}_f(\scru{\vec{x}}{gc}, k,a)$ of function $f(x,k)$ at the x-nodes with the coordinates $\scru{u}{gc}_i=a \scru{x}{gc}_i$ (\ProcName{CalcVectorVfGC}).}
	\label{alg-CalcVectorVfGC}
	\begin{algorithmic}[1]
		\Require Pointer to function $\PtrToFun{f}: \Real\times \Real \to \Real$, positive integer $n$, $n$-dim root vector $\scru{\vec{x}}{gc}$, positive $a$, and variable $k$. 
		\Ensure The vector $\scru{\vec{v}}{gc}_f(k,a)$.
		\Function{CalcVectorVfGC}{\PtrToFun{f}, $n$, $\scru{\vec{x}}{gc}$, $k$, $a$}
		\For{$i\in \set{1, \cdots, n}$}
		\State $\scru{v}{gc}_i \gets \PtrToFun{f}(a\cdot \scru{x}{gc}_i, k)$ ;
		\EndFor
		\State \Return $\scru{\vec{v}}{gc}$;
		\EndFunction
	\end{algorithmic}
\end{breakablealgorithm}

\begin{breakablealgorithm}
	\caption{Calculate the CEI-1 with the Gauss-Chebyshev method (\ProcName{Kintgc})}
	\label{alg-Kintgc}
	\begin{algorithmic}[1]
		\Require Variable $k\in [0,1)$, positive integer $n$ with default value $n = 50$
		\Ensure The value of CEI-1
		\Function{Kintgc}{$k, n$}
		\State $\scru{\vec{x}}{gc}\gets \ProcName{CalcRootsChebyshev}(n)$;
		\State $\scru{\vec{w}}{gc}\gets \ProcName{CalcWeightsGC}(\scru{\vec{x}}{gc},n)$;
		\State $\PtrToFun{f}\gets \ProcName{KerGC}$;
		\State $a\gets 1.0$;
		\State $K \gets \ProcName{IntegralSolverGC}(k, \PtrToFun{f}, n, \scru{\vec{x}}{gc}, \scru{\vec{w}}{gc}, a)$;
		\State \Return $K$;
		\EndFunction
	\end{algorithmic}
\end{breakablealgorithm}

\subsection{Gauss-Legendre Algorithm}

The key idea of the Gauss-Legendre method is the same as that of the Gauss-Chebyshev. However, it is more complex due to the following reasons:
\begin{itemize}
	\item it is necessary for us to compute the Legendre polynomials $P_n(x)$ and their derivatives $P'_n(x)$ automatically;
	\item the root vector $\scru{\vec{x}}{gl} = \trsp{\left[\scru{x}{gl}_1, \cdots,\scru{x}{gl}_n\right]}$ of $P_n(x)$ must be solved iteratively with the Newton's method;
	\item the weight vector $\scru{\vec{w}}{gl}=\trsp{\left[\scru{w}{gl}_1, \cdots, \scru{w}{gl}_n\right]}$ must be solved with the help of $\scru{\vec{x}}{gl}$, $P_n(x)$ and $P'_n(x)$;
	\item the Newton's method for solving roots of nonlinear equation should be designed.
\end{itemize}
The strategy for design algorithm for CEI-1 with Gauss-Legendre method is still top-down. In other words, we first design an algorithm to the general problem of computing the integral $\displaystyle{\int^b_a f(x,k)\dif x} = \braket{\scru{\vec{w}}{gl}}{\scru{\vec{v}}{gl}_f(\scru{\vec{x}}{gc}, k,a,b)} + \mathcal{E}_n(f)$ with 
inner product, then we calculate the vectors $\scru{\vec{x}}{gl}$ and $\scru{\vec{w}}{gl}$ with orthogonal polynomials and Newton's method for solving roots. 

There is also an alternative way to get the root vector $\scru{\vec{x}}{gl}$  and $\scru{\vec{w}}{gl}$ by reading some handbooks about special functions and/or numerical analysis. However, the tables obtained are just for some fixed significant digits and fixed order $n$ such that $n \le 20$ for $P_n(x)$. It is necessary for us to point out that a smart way born in 1900s is not always a wise way in 2020s. Moreover, the tables given in references can be used to check whether the output of the implementation of our algorithms is correct or not. It is a good idea that we respect the tradition, stand on the tradition, but also try to go beyond the tradition. 

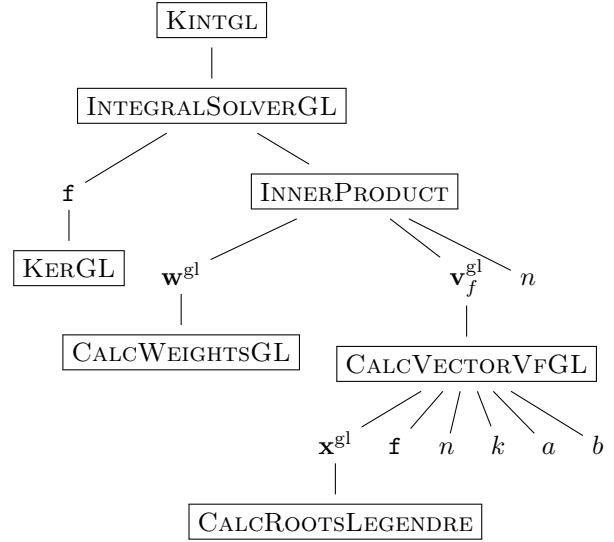
\begin{figure}[h] 
	\centering
	\begin{forest}
		[\fbox{\ProcName{Kintgl}}
		[\fbox{\ProcName{IntegralSolverGL}}
		[\PtrToFun{f}[\fbox{\ProcName{KerGL}}]]
		[\fbox{\ProcName{InnerProduct}}
		[$\scru{\vec{w}}{gl}$[\fbox{\ProcName{CalcWeightsGL}}]]
		[$\scru{\vec{v}}{gl}_f$
		[\fbox{\ProcName{CalcVectorVfGL}}
		[$\scru{\vec{x}}{gl}$
		[\fbox{\ProcName{CalcRootsLegendre}}]	        
		]
		[\PtrToFun{f}]	
		[$n$]
		[$k$]
		[$a$]
		[$b$]
		]	
		]
		[$n$]
		] 
		]	
		] 
	\end{forest}
	\caption{Procedures involved in calculating the complete elliptic integral CEI-1 with Gauss-Legendre method.}
	\label{fig-Kgl-procedures}
\end{figure}

The complete design for the algorithm of Gauss-Legendre for CEI-1 can be divided into four sub-tasks. We now discuss them one-by-one with a strategy of top-down. \Fig 
\ref{fig-Kgl-procedures} demonstrates the structure of  procedures in the Gauss-Legendre method. 

\subsection{Compute CEI-1 with the Gauss-Legendre based on the roots, weights and kernel function}.

Firstly, recognize the computation of CEI-1 as a special problem of computing the integral $\scrd{F}{L}(k)$ with numerical integration via  \eqref{eq-f-gl} and  \eqref{eq-K-gl}. The procedure \ProcName{IntegralSolverGL} in \Algr \ref{alg-IntegralSolverGL} is used to compute $\scrd{F}{L}(k)$ based on the
procedure \ProcName{InnerProduct} in \Algr \ref{alg-InnerProduct} and the procedure \ProcName{KerGL} in \Algr \ref{alg-KerGL} to compute the function $\scrd{f}{gl}(t,k) = 1/\sqrt{1-k^2\sin^2t}$. 

\begin{breakablealgorithm}
	\caption{Calcualte the parametric integral 
		$$\displaystyle 
		\scrd{F}{L}(y) =\int^b_a f(x,k)\dif x = \braket{\scru{\vec{w}}{gl}}{\scru{\vec{v}}{gl}_f(\scru{\vec{x}}{gl},k,a,b)} + \mathcal{E}_n(f)
		$$
		via the $n$-nodes Gauss-Legendre numerical integration method (\ProcName{IntegralSolverGL}).}
	\label{alg-IntegralSolverGL}
	\begin{algorithmic}[1]
		\Require Variable $k$, pointer to function $\PtrToFun{f}: \Real \times \Real \to \Real$, positive integral $n$, $n$-dim root vector $\scru{\vec{x}}{gl}$, $n$-dim weight vector $\scru{\vec{w}}{gl}$, variable $a$ and variable $b$.
		\Ensure The value  of $\scrd{F}{L}(y)$.
		\Function{IntegralSolverGL}{$k$, \PtrToFun{f}, $n$, $\scru{\vec{x}}{gl}$, $\scru{\vec{w}}{gl}$, $a$, $b$}
		\State $\scru{\vec{v}}{gl}\gets \ProcName{CalcVectorVfGL}( 
		\PtrToFun{f}, n, \scru{\vec{x}}{gl}, k, a, b)$;
		\State $\cpvar{j}\gets \ProcName{InnerProduct}(\scru{\vec{w}}{gl}, \scru{\vec{v}}{gl}, n)$;
		\State \Return \cpvar{j};
		\EndFunction
	\end{algorithmic}
\end{breakablealgorithm}

\begin{breakablealgorithm}
	\caption{Compute the function $\scrd{f}{gl}(t,k)$ in the Gauss-Legendre method for calculating the value CEI-1 (\ProcName{KerGL})}
	\label{alg-KerGL}
	\begin{algorithmic}[1]
		\Require Variables $t\in [-1,1], k\in [0,1)$
		\Ensure The value of $\scrd{f}{gl}(t,k)$
		\Function{KerGL}{$t,k$}
		\State $u \gets k\cdot \sin t$;
		\State $v\gets 1/\sqrt{1-u^2}$;
		\State \Return $v$;
		\EndFunction
	\end{algorithmic}
\end{breakablealgorithm}

\begin{breakablealgorithm}
	\caption{Calculate the vector $\scru{\vec{v}}{gl}_f(\scru{\vec{x}}{gl},k,a,b)$ of function $f(x,k)$ at the x-nodes with the coordinates $\scru{u}{gl}_i=\left[b+a + (b-a)\scru{x}{gc}_i\right]/2$ (\ProcName{CalcVectorVfGL}).}
	\label{alg-CalcVectorVfGL}
	\begin{algorithmic}[1]
		\Require The pointer to function $\PtrToFun{f}: \Real\times \Real \to \Real$, the order $n$, the $n$-dim root vector $\scru{\vec{x}}{gl}$, the variable $a$, the variable $b$, and the variable $k$.
		\Ensure Vector $\scru{\vec{v}}{gc}_f(\scru{\vec{x}}{gl},k,a,b)$.
		\Function{CalcVectorVfGL}{\PtrToFun{f}, $n$, $\scru{\vec{x}}{gc}$, $k$, $a$, $b$}
		\For{$i\in \seq{1, \cdots, n}$} \quad //$\scru{v}{gl}_{f,i}(\scru{\vec{x}}{gl},k,a,b)$;
		\State $\scru{u}{gl} \gets (b+a + (b-a)\scru{x}{gc}_i)/2$;
		\State $\scru{v}{gl}_i \gets \PtrToFun{f}(\scru{u}{gl}, k)\cdot \frac{b-a}{2}$; 
		\EndFor
		\State \Return $\scru{\vec{v}}{gl}$;
		\EndFunction
	\end{algorithmic}
\end{breakablealgorithm}

\begin{breakablealgorithm} 
	\caption{Calculate the CEI-1 with the Gauss-Legendre method (\ProcName{Kintgl})}
	\label{alg-Kintgl}
	\begin{algorithmic}[1]
		\Require The variable $k\in [0,1)$, the positive integer $n$ with default value $n = 50$
		\Ensure The value of CEI-1
		\Function{Kintgl}{$k, n$}
		\State $\scru{\vec{x}}{gl}\gets \ProcName{CalcRootsLegendre}(n)$;
		\State $\scru{\vec{w}}{gl}\gets \ProcName{CalcWeightsGL}(\scru{\vec{x}}{gl},n)$;
		\State $\PtrToFun{f}\gets \ProcName{KerGL}$;
		\State $a\gets 0$;
		\State $b\gets \pi/2$;
		\State $K \gets \ProcName{IntegralSolverGL}(k, f, n, \scru{\vec{x}}{gl},
		\scru{\vec{w}}{gl}, a, b)$;
		\State \Return $K$;
		\EndFunction
	\end{algorithmic}
\end{breakablealgorithm}

\subsection{Generating the Weights Specified by the Legendre polynomial and its Derivative}

\begin{figure}[h]
	\centering
	\begin{forest}
		[\fbox{\ProcName{CalcWeightsGL}}
		[$\scru{\vec{x}}{gc}$[\fbox{\ProcName{CalcRootsLegendre}}]]
		[$P'_n(x)$[\fbox{\ProcName{DerivPolynomLegendre}}]]
		]
	\end{forest}
	\caption{Procedures involved in computing the $n$ weights for the numerical integration with Gauss-Legendre method.} \label{fig-Kintgl-calc-weights}
\end{figure}
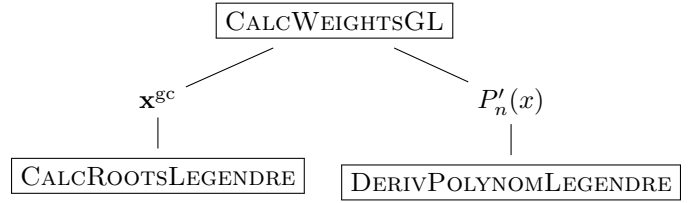

The weight vector $\scru{\vec{w}}{gl}$ for the Gauss-Legendre numerical integration can be generated with the help of the root vector $\scru{\vec{x}}{gl}$. 	According to \eqref{eq-weight-Legendre}, we design a procedure \ProcName{CalcWeightsGL} in \Algr \ref{alg-CalcWeightsGL} to calculate $\scru{\vec{w}}{gl}$ by calling the procedure \ProcName{DerivPolynomLegendre} in \Algr \ref{alg-DerivPolynomLegendre}. The details for the procedure \ProcName{CalcWeightsGL} are given in \Algr~\ref{alg-CalcWeightsGL}, which is a direct translation of \eqref{eq-weight-Legendre} from mathematics to pseudo-code. The structure of procedures for calculating the weights is shown in \Fig \ref{fig-Kintgl-calc-weights}.

\begin{breakablealgorithm}
	\caption{Calculate the weight vector $\scru{\vec{w}}{gl}$ for the abstract function (parametric integral with parameter $k$)
		$\scrd{F}{L}(k) = \displaystyle{\int^b_a f(x,k)\dif x} = \braket{\scru{\vec{w}}{gl}}{\scru{\vec{v}}{gl}_f(\scru{\vec{x}}{gl},k,a,b)} + \mathcal{E}_n(f)$ with the $n$-nodes Gauss-Legendre numerical integration method (\ProcName{CalcWeightsGL}).}
	\label{alg-CalcWeightsGL}
	\begin{algorithmic}[1]
		\Require Positive integer $n$ for the number of weights, root vector 
		$\scru{\vec{x}}{gl}$ of $P_n(x)$.
		\Ensure Weight vector $\scru{\vec{w}}{gl}$ for the Gauss-Legendre numerical integration.
		\Function{CalcWeightsGL}{$\scru{\vec{x}}{gl}, n$}
		\For{$i\in \seq{1, \cdots, n}$}
		\State $u \gets \ProcName{DerivPolynomLegendre}(\scru{x}{gl}_i,n)$;
		\State $\scru{w}{gl}_i \gets \cfrac{2}{[1-(\scru{x}{gl}_i)^2]u^2}$;
		\EndFor
		\State \Return $\scru{\vec{w}}{gl}$;
		\EndFunction
	\end{algorithmic}
\end{breakablealgorithm}

\subsection{Solving the $n$ Roots of the Legendre Polynomial $P_n(x)$}

\begin{figure}[htbp]
\centering
\includegraphics[width=0.5\textwidth]{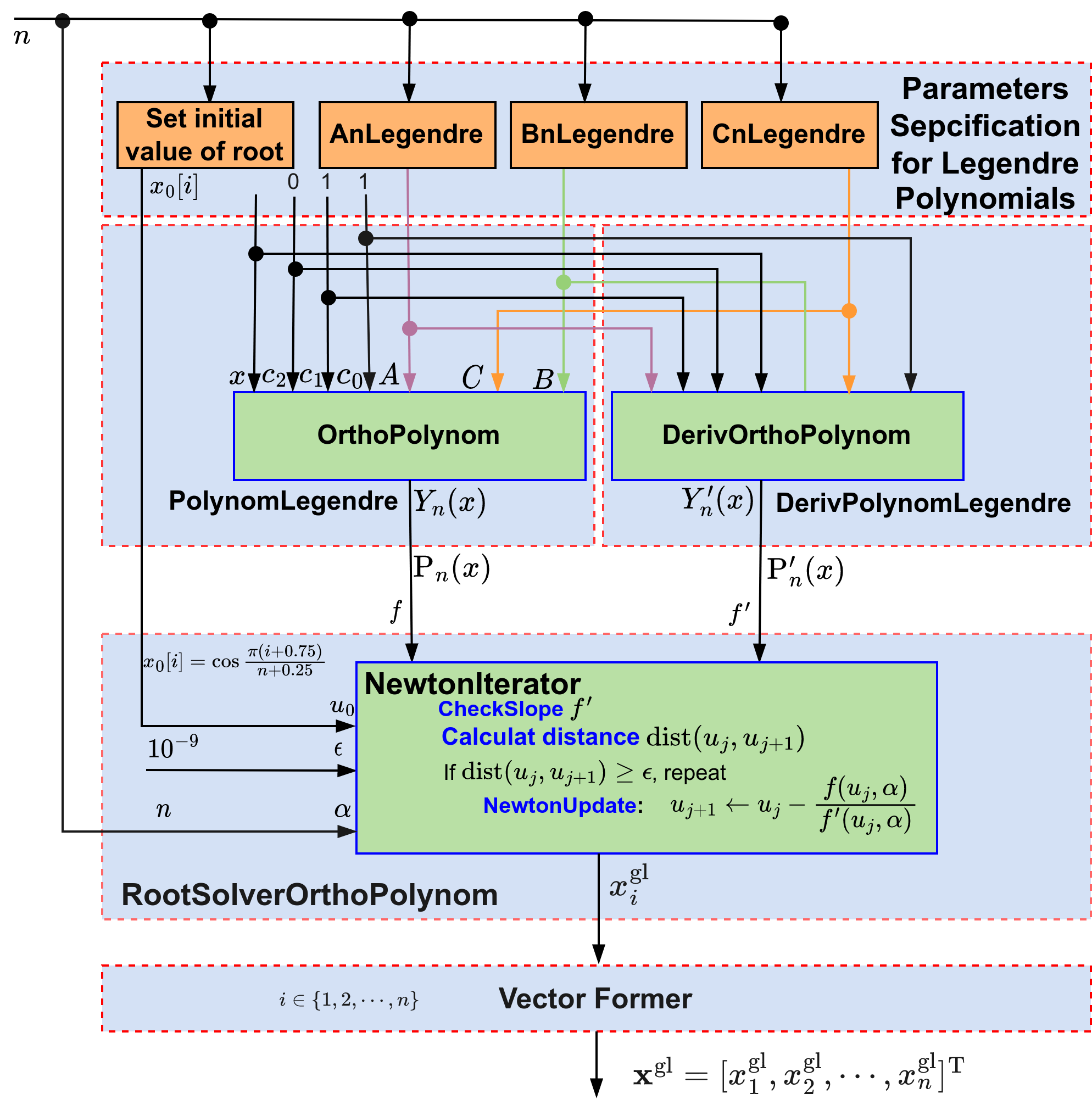} 
\caption{Solving the root vector $\scru{\vec{x}}{gl}$ for the Legendre polynomial $P_n(x)$ such that  $P_n(\scru{x}{gl}_{j})=0$.}
\label{fig-NewtonSolver}
\end{figure}

\Fig \ref{fig-NewtonSolver} demonstrates the system of solving the root vector of $P_n(x)$.  Usually in the books about special mathematical functions \cite{CourantHilbert-vol1,SJZhang1996,HandbookMPES2011}, the different kinds of orthogonal polynomials, such as the Legendre polynomials $P_n(x)$, Chebyshev polynomials $T_n(x)$ of the first or second kind, Laguerre polynomials $L_n(x)$ and Hermite polynomials $H_n(x)$,  are treated separately, which hides the general iterative formula \eqref{eq-ortho}. Once the $i$-th root $\scru{x}{gl}_i$ is obtained, the corresponding weight $\scru{w}{gl}_i$ can be determined by \eqref{eq-weight-Legendre}.

\begin{figure}[h] 
	\centering
	\begin{forest}
		[\fbox{\ProcName{CalcRootsLegendre}}
		[\fbox{\ProcName{RootSolverOrthoPolynom}}
		[\fbox{\ProcName{IsGoodEnough}}[\fbox{\ProcName{DistEuclid}}]]
		[\fbox{\ProcName{NewtonUpdate}}        		
		[\PtrToFun{Df}[\fbox{\ProcName{DerivPolynomLegendre}}]]
		[\PtrToFun{f}[\fbox{\ProcName{PolynomLegendre}}]]
		]
		]
		] 
	\end{forest}
	\caption{Procedures involved in calculating the $n$ roots of $P_n(x)$. }
	\label{fig-calc-roots-Legendre-procedures}
\end{figure}
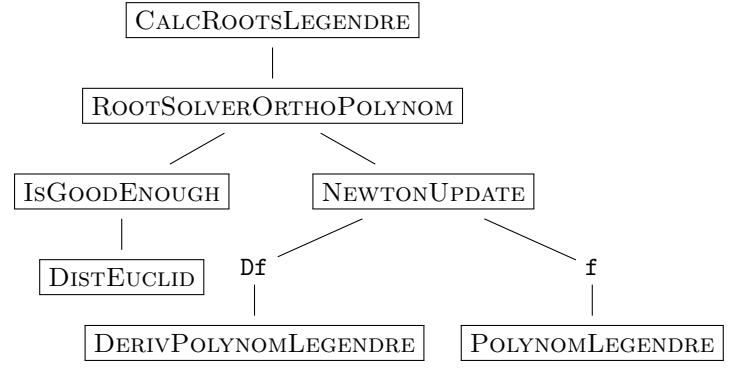

The steps for solving the $n$ roots of $P_n(x)$  lie in three issues:
\begin{itemize}
	\item designing the update procedure \ProcName{NewtonUpdate} in \Algr \ref{alg-NewtonUpdate}  by using $P_n(x)$ with the iteration formula in Newton's method;
	\item setting proper initial values $\scru{x}{gl}_i[0]$ for the $n$ roots of $P_n(x)$ according to  \eqref{eq-init-values-for-Pn} and call the procedure \ProcName{CalcRootsLegendre} in \Algr \ref{alg-CalcRootsLegendre};
	\item computing one root of $P_n(x)$ by the Newton's method encapsulated by the procedure \ProcName{RootSolverOrthoPolynom}  in \Algr \ref{alg-RootSolverOrthPolynom}.   
\end{itemize}
The structure of procedures for calculating the roots is illustrated in \Fig \ref{fig-calc-roots-Legendre-procedures}. 

\begin{breakablealgorithm}
	\caption{Calculte the root of orthogonal polynomial $Y_n(x)$ with an initial value $x_0$ via Newton's method (\ProcName{RootSolverOrthoPolynom}). }
	\label{alg-RootSolverOrthPolynom}
	\begin{algorithmic}[1]
		\Require Pointer to function  $\PtrToFun{f}: \Real \times \mathbb{N}\to \Real$ for $Y_n(x)$,  pointer to function $\PtrToFun{Df}: \Real \times \mathbb{N}\to \Real$ for the derivative function $Y'_n(x)$, initial value $x_0$, order $n$, and the precision $\varepsilon$.
		\Ensure The root $x$ of $Y_n(x) = 0$ iterated from $x_0$.
		\Function{RootSolverOrthoPolynom}{\PtrToFun{f}, \PtrToFun{Df}, $x_0$, $n$, $\varepsilon$} 
		\State $\cpvar{guess}\gets 0$ ;   
		\State $\cpvar{improve} \gets x_0$;
		\Repeat
		\State \cpvar{guess} $\gets$ \cpvar{improve};
		\State \cpvar{improve}   $\gets$ \ProcName{NewtonUpdate}(\PtrToFun{f}, \PtrToFun{Df}, \cpvar{guess}, $n$);
		\Until{$(\ProcName{DistEuclid}(\cpvar{guess}, \cpvar{improve},1)\ge \varepsilon)$}
		\State \Return \cpvar{improve};
		\EndFunction
	\end{algorithmic}
\end{breakablealgorithm} 

\begin{breakablealgorithm}
	\caption{Updating function $\scrd{\mathcal{A}}{Newton}$ for solving a root of the orthogonal polynomial $Y_n(x)$ with the Newton method  (\ProcName{NewtonUpdate}).}
	\label{alg-NewtonUpdate}
	\begin{algorithmic}[1]
		\Require Two pointers to function:  \PtrToFun{f} and \PtrToFun{Df} for $Y_n(x)$ and its derivative $Y'_n(x)$ respectively, variable $x$ and order $n\in \mathbb{N}$.
		\Ensure Updated value of the approximation for the root of interest
		\Function{NewtonUpdate}{\PtrToFun{f}, \PtrToFun{Df}, $x$, $n$}
		\State $\cpvar{slope}\gets \PtrToFun{Df}(x,n)$; // Compute  $Y'_n(x)$
		\State \ProcName{CheckSlope}(\cpvar{slope});  // Check $Y'_n(x)$;
		\State $\scrd{x}{new}  \gets x - \cfrac{\PtrToFun{f}(x, n)}{\cpvar{slope}}$; 
		// $x_{i+1} = x_i - \cfrac{Y_n(x_i)}{Y'_n(x_i)}$
		\State \Return $\scrd{x}{new} $;
		\EndFunction
	\end{algorithmic}
\end{breakablealgorithm}

\begin{breakablealgorithm} 
	\caption{Check whether the slope/derivative at point $x$ is zero or not (\ProcName{CheckSlope}).}
	\label{alg-CheckSlope}
	\begin{algorithmic}[1]
		\Require The slope (value of derivative function) $\cpvar{slope}$
		\Ensure Void or exception information if necessary
		\Function{CheckSlope}{\cpvar{slope}}
		\State $\delta \gets 10^{-10}$;
		\If{($\abs{\cpvar{slope}}< \delta$)}
		\State  \cpvar{Print Exception Information!}
		\State  \textbf{Exit};
		\EndIf
		\EndFunction
	\end{algorithmic}
\end{breakablealgorithm}

\begin{breakablealgorithm}
	\caption{Calculate the $n$-dim root vector  $\scru{\vec{x}}{gl}=\trsp{[\scru{x}{gl}_1, \cdots, \scru{x}{gl}_n]}$ of the $n$--th order Legendre polynomial $P_n(x)$ such that
		$P_n(\scru{x}{gl}_i) = 0$ for $1\le i\le n$ (\ProcName{CalcRootsLegendre}).}
	\label{alg-CalcRootsLegendre}
	\begin{algorithmic}[1]
		\Require The positive integer $n$
		\Ensure The root vector $\scru{\vec{x}}{gl}$ of $P_n(x)$
		\Function{CalcRootsLegendre}{$n$}
		\State Allocate memory for $\scru{\vec{x}}{gl}=\trsp{[\scru{x}{gl}_1, \cdots, \scru{x}{gl}_n]}$;
		\State $\varepsilon \gets 10^{-9}$; //set the precision for root
		\State $\PtrToFun{f} \gets \ProcName{PolynomLegendre}$;  // $f (x)= P_n(x)$
		\State $\PtrToFun{Df} \gets \ProcName{DerivPolynomLegendre}$; // $f'=P'_n(x)$ 
		\For{$i\in \seq{1,2,\cdots,n}$}
		\State $x_0 \gets \cos \cfrac{\pi(i+0.75)}{n+0.25}$; //initial value; 
		\State $\scru{x}{gl}_i\gets$ \ProcName{RootSolverOrthoPolynom}(\PtrToFun{f}, \PtrToFun{Df}, $x_0$, $n$, $\varepsilon$);
		\EndFor
		\State \Return $\scru{\vec{x}}{gl}$;
		\EndFunction
	\end{algorithmic}
\end{breakablealgorithm}

\subsection{Generating the Legendre polynomials and their derivatives}

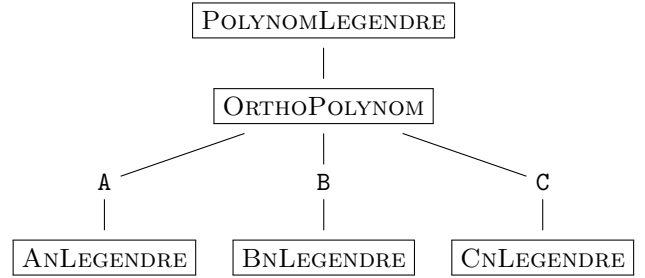
\begin{figure}[htp] 
	\centering
	\begin{forest}
		[\fbox{\ProcName{PolynomLegendre}}
		[\fbox{\ProcName{OrthoPolynom}}
		[\PtrToFun{A}[\fbox{\ProcName{AnLegendre}}]]
		[\PtrToFun{B}[\fbox{\ProcName{BnLegendre}}]]
		[\PtrToFun{C}[\fbox{\ProcName{CnLegendre}}]]
		]
		] 
	\end{forest}
	\caption{Procedures for calculating $P_n(x)$ with a top-down strategy: \PtrToFun{A}, \PtrToFun{B} and \PtrToFun{C} are pointers to function with the assignments $\PtrToFun{A}\gets \ProcName{AnLegendre}$, $\PtrToFun{B}\gets \ProcName{BnLegendre}$ and $\PtrToFun{C}\gets \ProcName{CnLegendre}$ respectively. }
	\label{fig-procs-Pn}
\end{figure}

\begin{figure}[htp] 
	\centering
	\begin{forest}
		[\fbox{\ProcName{DerivPolynomLegendre}}
		[\fbox{\ProcName{DerivOrthoPolynom}}
		[\PtrToFun{A}[\fbox{\ProcName{AnLegendre}}]]
		[\PtrToFun{B}[\fbox{\ProcName{BnLegendre}}]]
		[\PtrToFun{C}[\fbox{\ProcName{CnLegendre}}]]
		]
		] 
	\end{forest}
	\caption{Procedures for calculating $P'_n(x)$ with a top-down strategy: \PtrToFun{A}, \PtrToFun{B} and \PtrToFun{C} are pointers to function with the assignments $\PtrToFun{A}\gets \ProcName{AnLegendre}$, $\PtrToFun{B}\gets \ProcName{BnLegendre}$ and $\PtrToFun{C}\gets \ProcName{CnLegendre}$ respectively. }
	\label{fig-procs-Deriv-Pn}
\end{figure}
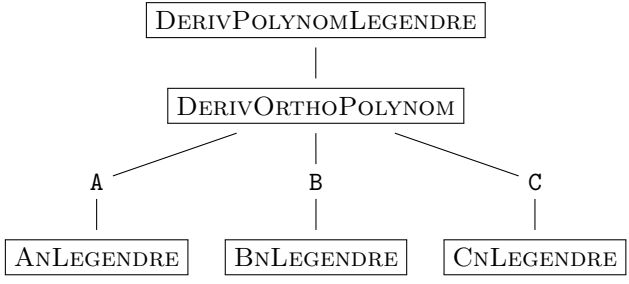

The general form of four typical different orthogonal polynomials (such as Legendre polynomials $P_n(x)$, Chebyshev polynomials $T_n(x)$, Laguerre polynomials $L_n(x)$ and Hermite polynomials $H_n(x)$), denoted by $Y_n(x)$,  can be specified by the following iterative formula
\begin{equation}
	\left\{
	\begin{split}
		Y_n(x) & = (A_nx + B_n)Y_{n-1}(x) - C_n Y_{n-2}(x), \quad  n\ge 2	 \\
		Y_0(x) & = c_0 \\
		Y_1(x) & = c_1 x + c_2
	\end{split}
	\right.
\end{equation}
where $c_0, c_1, c_2$ are constants, $A_n, B_n, C_n$ are functions depending on $n$. 
Taking the derivative on the two sides, we immediately have
\begin{equation} \label{eq-ortho}
	\left\{
	\begin{split}
		Y'_n(x) &= A_n Y_{n-1}(x) + (A_nx + B_n)Y'_{n-1}(x)
		- C_n Y'_{n-2}(x) \\
		Y'_0(x) & = 0 \\
		Y'_1(x) & = c_1
	\end{split}
	\right.
\end{equation}
where $n\ge 2$. 
Particularly, let
\begin{equation}
	\left\{
	\begin{split}
		\scru{A}{le}_n &= \frac{2n-1}{n}, \scru{B}{le}_n = 0, \scru{C}{le}_n = \frac{n-1}{n} \\
		\scru{c}{le}_0 &= 1,\scru{c}{le}_1 = 1, \scru{c}{le}_2 = 0
	\end{split}
	\right.
\end{equation}
then for $n\ge 2$ we can deduce that 
\begin{equation}
	\left\{
	\begin{split}
		P_n(x) & = \scru{A}{le}_nxP_{n-1}(x) - \scru{C}{le}_n P_{n-2}(x)	 \\
		P_0(x) & = 1 \\
		P_1(x) & = x
	\end{split}
	\right.
\end{equation}
and 
\begin{equation}
	\left\{
	\begin{split}
		P'_n(x) &= \scru{A}{le}_n P_{n-1}(x) + \scru{A}{le}_nxP'_{n-1}(x)
		- \scru{C}{le}_n P'_{n-2}(x) \\
		P'_0(x) & = 0 \\
		P'_1(x) & = 1
	\end{split}
	\right.
\end{equation}
This implies that we can compute the $n$-th order Legendre polynomial $P_n(x)$ and its derivative $P'_n(x)$ iteratively with a  strategy of top-down:
\begin{itemize}
	\item firstly, we design general procedures \ProcName{OrthoPolynom} in \Algr \ref{alg-OrthoPolynom} and \ProcName{DerivOrthoPolynom} in \Algr \ref{alg-DerivOrthoPolynom} for generating $Y_n(x)$ and $Y'_n(x)$ respectively with the iterative formulas directly;
	\item secondly, we design specific procedures --- \ProcName{AnLegendre} in \Algr \ref{alg-AnLegendre}, \ProcName{BnLegendre} in \Algr \ref{alg-BnLegendre} and \ProcName{CnLegendre} in \Algr \ref{alg-CnLegendre} --- to calculate $A(n), B(n)$ and $C(n)$ respectively and use pointer to functions to compute the Legendre polynomials $P_n(x)$ and $P'_n(x)$.
\end{itemize} 
The procedures involved in computing the Legendre polynomial $P_n(x)$ and its derivative $P'_n(x)$ are illustrated in \Fig\ref{fig-procs-Pn} and  \Fig \ref{fig-procs-Deriv-Pn} respectively.

\begin{breakablealgorithm} 
	\caption{Compute the orthogonal polynomials $Y_n(x)$ according to the iterative formula. (\ProcName{OrthoPolynom})}
	\label{alg-OrthoPolynom}
	\begin{algorithmic}[1]
		\Require The variable $x$, order $n\in \mathbb{N}$, three constants$\set{c_0, c_1, c_2}$, three pointers to function  $\PtrToFun{A}, \PtrToFun{B},\PtrToFun{C}: \mathbb{N}\to \Real$.
		\Ensure The value of $Y_n(x)$
		\Function{OrthoPolynom}{$x,n,c_0,c_1,c_2, \PtrToFun{A},\PtrToFun{B},\PtrToFun{C})$}
		\State $Y_0 \gets c_0$;   // $Y_0(x) = c_0$;
		\State $Y_1 \gets c_1 x + c_2$;  // $Y_1(x) = c_1 x +c_2$;
		\State $Y_n \gets 0$;  // Initialize $Y_n(x)$ with $0$;
		\If{$(n == 0)$}  \State $Y_n \gets Y_0$; \EndIf
		\If{$(n == 1)$}  \State $Y_n \gets Y_1$; \EndIf
		\For{$k\in \seq{2, \cdots, n}$}
		\State $Y_n \gets (\PtrToFun{A}(k) \cdot x + \PtrToFun{B}(k))\cdot Y_1 - \PtrToFun{C}(k)\cdot Y_0$;
		\State $Y_0\gets Y_1$; // Update $Y_{k-2}(x)$ with $Y_{k-1}(x)$;
		\State $Y_1\gets Y_n$; // Update $Y_{k-1}(x)$ with $Y_{k}(x)$;
		\EndFor
		\State \Return $Y_n$;
		\EndFunction
	\end{algorithmic}
\end{breakablealgorithm}

\begin{breakablealgorithm} 
	\caption{Compute the derivative $Y'_n(x)$ of the orthogonal polynomial $Y_n(x)$ according to the
		iterative formula (\ProcName{DerivOrthoPolynom}).}
	\label{alg-DerivOrthoPolynom}
	\begin{algorithmic}[1]
		\Require The variable $x$, the order $n\in \mathbb{N}$, the constants $c_0, c_1$ and $c_2$, the three pointers to function  \PtrToFun{A}, \PtrToFun{B} and \PtrToFun{C}.
		\Ensure The value of the derivative function $Y'_n(x)$
		\Function{DerivOrthoPolynom}{$x$,$n$,$c_0$,$c_1$,$c_2$,\PtrToFun{A},\PtrToFun{B},\PtrToFun{C}}
		\State $Y_0 \gets c_0$;  // $Y_0(x) = c_0$
		\State $Y_1 \gets c_1 x + c_2$;  // $Y_1(x) =c_1 x + c_2$
		\State $ Y'_0 \gets 0$; // $Y_0(x) = c_1 \Longrightarrow Y'_0(x)=0$
		\State $Y'_1 \gets c_1$; // $ Y_1(x) = c_1x+c_2 \Longrightarrow Y'_1(x)=c_1$
		\State $Y_n \gets 0$;  // Initialize $Y_n(x)$ with $0$
		\State $Y'_n \gets 0$;  // Initialize $Y'_n(x)$ with $0$;
		\If{$(n == 0)$}
		\State $Y_n \gets Y_0$;
		\State $Y'_n \gets Y'_0$;
		\EndIf
		\If{$(n == 1)$}
		\State $Y_n \gets Y_1$;
		\State $Y'_n \gets Y'_0$;
		\EndIf
		\For{$k\in \seq{2, \cdots, n}$}
		\State  $Y_n \gets (\PtrToFun{A}(k) \cdot x + \PtrToFun{B}(k))\cdot Y_1 - \PtrToFun{C}(k)\cdot Y_0$;
		\State  $Y'_n \gets \PtrToFun{A}(k)\cdot Y_1 + (\PtrToFun{A}(k) \cdot x + \PtrToFun{B}(k))\cdot Y'_1  - \PtrToFun{C}(k)\cdot Y'_0$;
		\State  $Y_0\gets Y_1$;  // Update $Y_{k-2}(x)$ with $Y_{k-1}(x)$;
		\State $Y_1\gets Y_n$;   // Update $Y_{k-1}(x)$ with $Y_k(x)$;
		\State $Y'_0\gets Y'_1$; // Update $Y'_{k-2}(x)$ with $Y'_{k-1}(x)$; 
		\State  $Y'_1\gets Y'_n$;// Update $Y'_{k-1}(x)$ with $Y'_k(x)$;
		\EndFor
		\State \Return $Y'_n$;
		\EndFunction
	\end{algorithmic}
\end{breakablealgorithm}

\begin{breakablealgorithm} 
	\caption{Compute the coefficient $A(n)$ in the iterative formula for Legendre polynomial (\ProcName{AnLegendre}).}
	\label{alg-AnLegendre}
	\begin{algorithmic}[1]
		\Require Order $n\in \mathbb{N}$
		\Ensure The value of $\scru{A}{le}_n$
		\Function{AnLegendre}{$n$}
		\State  $A\gets 2.0 -1.0/n$;
		\State \Return $A$;
		\EndFunction
	\end{algorithmic}
\end{breakablealgorithm}

\begin{breakablealgorithm} 
	\caption{Compute the coefficient $B(n)$ in the iterative formula for Legendre polynomial (\ProcName{BnLegendre}).}
	\label{alg-BnLegendre}
	\begin{algorithmic}[1]
		\Require Order $n\in \mathbb{N}$
		\Ensure The value of $\scru{B}{le}_n$
		\Function{BnLegendre}{$n$}
		\State \Return $0$;
		\EndFunction
	\end{algorithmic}
\end{breakablealgorithm}

\begin{breakablealgorithm} 
	\caption{Compute the coefficient $C(n)$ in the iterative formula for Legendre polynomial (\ProcName{CnLegendre}).}
	\label{alg-CnLegendre}
	\begin{algorithmic}[1]
		\Require Order $n\in \mathbb{N}$
		\Ensure The value of $\scru{C}{le}_n$
		\Function{CnLegendre}{$n$}
		\State $C\gets 1.0 -1.0/n$;
		\State \Return $C$;
		\EndFunction
	\end{algorithmic}
\end{breakablealgorithm}

\begin{breakablealgorithm}
	\caption{Compute the Legendre Polynomial $P_n(x)$ (\ProcName{PolynomLegendre}).}
	\label{alg-PolynomLegendre}
	\begin{algorithmic}[1]
		\Require The variable $x$, order $n\in \mathbb{N}$
		\Ensure The value of $P_n(x)$ 
		\Function{PolynomLegendre}{$x,n$}
		\State $c_0\gets 1, c_1\gets 1, c_2 \gets 0$;
		\State $\PtrToFun{A}\gets \ProcName{AnLegendre}$;
		\State $\PtrToFun{B}\gets \ProcName{BnLegendre}$;
		\State $\PtrToFun{C}\gets \ProcName{CnLegendre}$;
		\State $P_n \gets \ProcName{OrthoPolynom}(x,n,c_0,c_1,c_2,\PtrToFun{A},\PtrToFun{B},\PtrToFun{C})$;
		\State \Return $P_n$;
		\EndFunction
	\end{algorithmic}
\end{breakablealgorithm}

\begin{breakablealgorithm}
	\caption{Compute the derivative of Legendre polynomial (\ProcName{DerivPolynomLegendre}).}
	\label{alg-DerivPolynomLegendre}
	\begin{algorithmic}[1]
		\Require The variable $x$, order $n\in \mathbb{N}$
		\Ensure  The derivative of $P_n(x)$, i.e., the value of $P'_n(x)$ 
		\Function{DerivPolynomLegendre}{$x,n$}
		\State $c_0\gets 1$;
		\State $c_1\gets 1$;
		\State $c_2 \gets 0$;
		\State $\PtrToFun{A}\gets \ProcName{AnLegendre}$;
		\State $\PtrToFun{B}\gets \ProcName{BnLegendre}$;
		\State $\PtrToFun{C}\gets \ProcName{CnLegendre}$;
		\State $P'_n \gets \ProcName{DerivOrthoPolynom}(
		x,n,c_0,c_1,c_2,\PtrToFun{A},\PtrToFun{B},\PtrToFun{C})$;
		\State \Return $P'_n$;
		\EndFunction
	\end{algorithmic}
\end{breakablealgorithm}

It should be remarked that designing is a mixture of multiple factors in a general viewpoint.  Our discussion is limited in a narrow perspective which emphasizing algorithm design. For more deep thinking about design, please see Brooks' famous book on design\cite{Brooks2010}.

\section{Verification, Validation and Testing} \label{sec-vvt}

For the algorithms designed for computing the CEI-1, it is  necessary for us to construct some simple but useful test cases for verification, validation and testing.

The verification of the four methods for solving CEI-1 is simple if all of the modules involved have been verified and validated. In the Appendix \ref{sec-testcases-VVT}, we illustrate some typical test cases and omit the details for testing the modules for constructing the four numerical solvers for CEI-1 since it is  simple and straight forward. After that, what we need do is just to compare the values of $K(k)$ obtained by \Algr \ref{alg-Kseries}, \Algr \ref{alg-Kagm}, \Algr \ref{alg-Kintgc} and \Algr \ref{alg-Kintgl}. 

\begin{table}[htb]
\centering
\caption{A comparison of numerical solvers for $K(k)$ and  $\scrd{\mathcal{K}}{cs}(m)$ in MATLAB and Mathematica such that $m=k^2$}
\label{tab-comparison-K-solvers}
\begin{tabular}{cccc}
\hline
\makecell[c]{\textbf{Parameter} \\ $\quad$ \\ $k$ }  & \makecell[c]{\textbf{CEI-1} $K(k)$ \\ $\ProcName{Kseries},\ProcName{Kagm}$ \\ $\ProcName{Kintgc},\ProcName{Kintgl}$ } & \makecell[c]{\textbf{Parameter}\\  \\  $m$ } & \makecell[c]{\textbf{CEI-1} $\scrd{\mathcal{K}}{cs}(m)$ \\ \lstinline|ellipticK| \\ \lstinline|EllipticK|} \\
\hline
$0.00$  &    1.5707963    &    $0.00$     &  1.5707963  \\
$0.10$  &    1.5747456    &    $0.01$     &  1.5747456  \\
$0.20$  &    1.5868678    &    $0.04$     &  1.5868678  \\
$0.30$  &    1.6080486    &    $0.09$     &  1.6080486  \\
$0.40$  &    1.6399999    &    $0.16$     &  1.6399999  \\
$0.50$  &    1.6857504    &    $0.25$     &  1.6857504  \\ 
\hline
\end{tabular}
\end{table}

\Tab \ref{tab-comparison-K-solvers} shows the difference of numerical solvers. It is easy to find that
\begin{itemize}
\item Our four numerical algorithms lead to the same output and the correctness can be verified by the $\ProcName{Kseries}$ since it is based on the expansion of formula \eqref{eq-K-expansion} which is right definitely and can be handled manually.
\item The available funcitons \lstinline|ellipticK| in MATLAB and
\lstinline|EllipticK| in Mathematica have consistent numerical solutions, which implies that the algorithms adopted may be the same.
\end{itemize}
 It should be remarked that the built-in function $\scrd{\mathcal{K}}{cs}(m)$ in MATLAB and Mathematica are different from the $K(k)$ in our implementation although they are equivalent essentially.

\begin{table*}[h]
\centering
\caption{Validation of CEI-1 via Precision of Computation}
\label{tab-comp-precision}
\begin{tabular}{p{2.5cm}p{3cm}p{8cm}}
\hline
\textbf{Type of Solver}  & \makecell[l]{\textbf{Algorithm} \\ \textsc{Procedure}} &   \makecell[l]{\textbf{Default Value of Precision}: Relative or Absolute Error \\ Configuration of Parameters for Numerical Integration}  \\
\hline \hline
IS-CEI-1 & \makecell[l]{\Algr \ref{alg-Kseries} \\ $\ProcName{Kseries}(k,\varepsilon)$}   &   Relative error  $\varepsilon = 10^{-9}$ \\
\hline
AGM-CEI-1 & \makecell[l]{\Algr \ref{alg-AGM} \\ $\ProcName{Kagm}(k,\varepsilon)$}
&  Absolute error $\varepsilon = 10^{-9}$ \\
\hline
GC-CEI-1  & \makecell[l]{\Algr \ref{alg-Kintgc}  \\ $\ProcName{Kintgc}(k,n)$}    &  \makecell[l]{Absolute error: $\abs{\scru{\mathcal{E}}{gc}_n(\scrd{f}{gc},1,k)} < \varepsilon = 10^{-9} $  \\
 $(k,n) \in \set{(0.1, 4), (0.2, 5), (0.3, 6), (0.4, 7), (0.5, 8)}$ }   \\
\hline
GL-CEI-1  & \makecell[l]{\Algr \ref{alg-Kintgl} \\ $\ProcName{Kintgl}(k,n)$}    &  \makecell[l]{Absolute error: $\abs{\scru{\mathcal{E}}{gl}_n(\scrd{f}{gl},0,\pi/2,k)} < \varepsilon = 10^{-9} $   \\  $(k,n) = (0.1, 2)$ }   \\
\hline    
\end{tabular}
\end{table*}

We remark that the validation of the four methods for computing $K(k)$ should be noted because there are different parameters which control the precision of computation. As shown in \Tab \ref{tab-comp-precision}, it is not difficult to find the following features: 
\begin{itemize}
\item For the IS-CEI-1 solver, the meaning of relative precision $\varepsilon$ is clear, and the absolute error is also applicable. 
\item For the AGM-CEI-1 solver, the meaning of absolute error is clear, and the relative error could be adopted easily.  
\item For the GC-CEI-1 solver, $n$ is an implicit parameter for control the precision. Given the absolute error $\varepsilon$ and the parameter $k$, we should estimate the integer $\scru{n}{gc}_f$ mentioned in Appendix \ref{appendix-ErrEstGC} such that
\begin{equation} 
\abs{\scru{\mathcal{E}}{gc}_n(f,1,k)}<\varepsilon
\end{equation}
for $n > \scru{n}{gc}_f$. 
\item For the GL-CEI-1 solver, $n$ is also an implicit parameter for control the precision. Given the absolute error $\varepsilon$ and the parameter $k$, we should estimate  the integer $\scru{n}{gl}_f$ mentioned in Appendix \ref{appendix-ErrEstGL} such that 
\begin{equation} 
\abs{\scru{\mathcal{E}}{gl}_n(f, 0, \pi/2,k)}<\varepsilon
\end{equation}
for $n > \scru{n}{gl}_f$. 
\end{itemize}

Why do we use the implicit parameter $n$ instead of the 
explicit parameter $\varepsilon$ to control the precision of computation in the GC-CEI-1 and GL-CEI-1 methods? The answer lies in the difficulty of estimating the upper bounds of $\abs{\scrd{f}{gc}^{(2n)}(\xi,k)}$ and  $\abs{\scrd{f}{gl}^{(2n)}(\xi,k)}$ for $\xi \in (-1,1)$ and $k\in [0,\sqrt{2}/2]$. Consequently, it is difficult for us to design algorithms to estimate the  $\scru{n}{gl}_f$ and $\scru{n}{gc}_f$ with the precision $\varepsilon$ directly. This feature gives us a good example to demonstrate the advantages and disadvantages of the methods and algorithms developed.
It is valuable in the sense of education although it is a pity in the sense of engineering.

\section{Conclusions} \label{sec-conclusions}

The numerical schemes for computing the CEI-1 is essential for its applications. The available numerical solutions obtained by MATLAB and Mathematica are different due to the different of definitions for $K(k)$ and $\scrd{\mathcal{K}}{cs}(m)$, which implies that  we should pay attention to the verification, validation and testing stage  in research, development and STEM education in order to avoid potential errors.

Theoretical analysis shows that the correctness of the numerical solution can be verified by the analytic formula in the form of infinite series. As alternative choices, we proposed four feasible numerical schemes and designed  the \ProcName{Kseries}, \ProcName{Kagm}, \ProcName{Kintgc} and \ProcName{Kintgl} solvers for CEI-1 to obtain consistent numerical results based on infinite series method, AGM method, Gauss-Chebyshev \& Gauss-Legendre numerical integrations. The key algorithms developed have been verified, validated and tested.

Although our emphasis is put on the numerical algorithms for the CEI-1,
the methods adopted can also be  used for the complete elliptic integrals of the second and third kinds. Furthermore, the key modules for the four solvers are designed for  more general problems so as to satisfy various applications as much as possible. Particularly, the iterative algorithms for generating orthogonal polynomials with iterative formula $Y_n(x) = [A(n)x + B(n)]Y_{n-1}(x) - C(n) Y_{n-2}(x)$  and their derivatives can be taken to compute the  Legendre, Chebyshev, Laguerre and Hermite polynomials which have wide applications in mathematics, physics and engineering.

\appendix

\section{Error Estimation of Gaussian Numerical Integration} \label{appendix-ErrEst}

\subsection{Gauss-Chebyshev Formula and Approximation Error} 
\label{appendix-ErrEstGC}

The numerical integration was discussed by Zhang \cite{SJZhang1996} and 
Atkinson \cite{Atkinson2009}.
\begin{thm}
For $f\in C^{2n}[-1,1]$, the integral $\displaystyle{\int^1_{-1} \frac{f(x)}{\sqrt{1-x^2}}\dif x }$ can be computed by the Gauss-Chebyshev formula
\begin{equation*}
\int^1_{-1} \frac{f(x)}{\sqrt{1-x^2}}\dif x = \sum_{i=1}^n \scru{w}{gc}_i f(\scru{x}{gc}_i) + \scru{\mathcal{E}}{gc}_n(f).
\end{equation*}
where 
\begin{equation}
\scru{\mathcal{E}}{gc}_n(f) = \frac{2\pi}{2^{2n}(2n)!}f^{(2n)}(\eta), \quad \eta\in (-1,1)
\end{equation}
is the approximation error.
\end{thm}

By changing the variable $x$ with $t = x/a$, we can obtain
\begin{equation} \label{eq-integral-gc}
\int^a_{-a} \frac{f(x, k)}{\sqrt{a^2 - x^2}}\dif x 
= \int^1_{-1} \frac{f(at, k)}{\sqrt{1-t^2}}\dif t.
\end{equation}
According to \eqref{eq-integral-gc} and the $k$-th order derivative formula 
\begin{equation} \label{eq-compos-deriv}
\left[\fracode{}{x}\right]^k f(\alpha x + \beta) = \alpha^k f^{(k)}(\alpha x + \beta),
\end{equation}
we can deduce the following corollary.

\begin{cor}
For $f(\cdot, k)\in C^{2n}[-1,1]$ and fixed parameter $k$, the integral $\displaystyle{\int^a_{-a} \frac{f(x,k)}{\sqrt{a^2-x^2}}\dif x}$ can be computed with the extended Gauss-Chebyshev formula
\begin{equation} \label{eq-F-Chebyshev-sum}
\int^a_{-a} \frac{f(x,k)}{\sqrt{a^2-x^2}}\dif x = \sum^{n}_{i=1} \scru{w}{gc}_i f(a\scru{x}{gc}_i, k) +  \scru{\mathcal{E}}{gl}_n(f,a,k)
\end{equation}
where  $\xi \in (-a, a)$ and 
\begin{equation}
\scru{\mathcal{E}}{gc}_n(f,a,k) = \frac{2\pi\cdot a^{2n}}{2^{2n}(2n)!}f^{(2n)}(\xi,k)
\end{equation}
is the approximation error.
\end{cor}

Theoretically, if there is an upper bound $M^f_n(k)$ for $f^{(2n)}(\xi,k)$ such that $\abs{f^{(2n)}(\xi,k)}\le M_n(k)$ where $\xi \in (-a,a)$, then we can find an integer 
\begin{equation}
\scru{n}{gc}_f = \scru{n}{gc}_f(k,\epsilon,a) =  \arg \min_{n\in \mathbb{N}}\set{n: \frac{2\pi a^{2n}}{2^{2n}(2n)!}M^f_n(k) < \epsilon}
\end{equation}
such that $\scru{\mathcal{E}}{gc}_n(f,a,k) < \epsilon$ when $n > \scru{n}{gc}_f$.

For illustration, if we take $f(x) = \cos x$, then $f^{(2n)}(x) = (-1)^n\cos x, f^{(2n)}(kx) = k^{2n}(-1)^n\cos kx$. Therefore, the  approximation error of integral $\displaystyle{\int^1_{-1}\frac{\cos kx}{\sqrt{1-x^2}}\dif x} = J_0(\abs{k}), \quad k\in \Real$  with the Gauss-Chebyshev formula will be 
\begin{align*}
\scru{\mathcal{E}}{gc}_n(\cos, 1,k) 
&= \frac{2\pi}{2^{2n}(2n)!}k^{2n}f^{(2n)}(\xi,k), \quad \xi \in (-1,1) \\
&= \frac{2\pi}{2^{2n}(2n)!}k^{2n}(-1)^n\cos(k \xi)
\end{align*} 
For the given $\epsilon > 0$, the upper bound of the absolute approximate error $\abs{\scru{\mathcal{E}}{gc}_n(\cos, 1,k)}< \epsilon$ means that $M^{\cos}_n(k) = k^{2n}$ and 
\[
\frac{2\pi}{2^{2n}(2n)!}k^{2n} < \epsilon 
\]
since $\abs{(-1)^n\cos(k \xi)}\le 1$.  For the fixed $k$, we can set the parameter $n$ as 
\[
\scru{n}{gc}_{\cos} = \arg \min_{i\in \mathbb{N}} \set{i: \frac{2\pi\cdot k^{2i}}{2^{2i}(2i)!}<\epsilon}
\]
Similarly, for $f(x) = \me^{x}$, we can deduce that 
\[
\scru{n}{gc}_{\me} = \arg \min_{i\in \mathbb{N}} \set{i: \frac{2\pi\cdot k^{2i}\me^k}{2^{2i}(2i)!}<\epsilon}
\]

For $\epsilon = 10^{-9}$, we have the $(k,n)$ pairs demonstrated in \Tab \ref{tab-n-gc-gl}

\begin{table}[H]
\begin{center}
\caption{Comparison of $\scru{n}{gc}_{\cos}$ and $\scru{n}{gc}_{\me}$}
\label{tab-n-gc-gl}
\begin{tabular}{|c|ccccccccc|}
\hline
$k$                  & 0.1  & 0.5  & 1.5 & 5 & 10 & 50 & 100 & 200 & 400\\
\hline
$\scru{n}{gc}_{\cos}$ &  4  &  6    &  7   &  8 &  10  &  14  &  18   & 22    &  28  \\
$\scru{n}{gc}_{\me}$ &  5  &   6   &   7  &  9  &  12  &  26  &  42   &  69   &  117  \\
\hline
\end{tabular}
\end{center}
\end{table}

In general, we can find the $\scru{n}{gc}_f$ such that the absolute approximation error is less than $\epsilon$ for the fixed parameter $k$ if the upper bound of $\abs{f^{(2n)}(\xi,k)}$ can be estimated properly.

\subsection{Gauss-Legendre Formula and Approximation Error}
\label{appendix-ErrEstGL}

\begin{thm} \label{thm-approx-error-GL}
For $f\in C^{2n}[-1,1]$, the integral $\displaystyle{\int^1_{-1}f(x)\dif x} $ can be computed with the Gauss-Legendre formula
\begin{equation*}
\int^1_{-1} f(x)\dif x = \sum^n_{i=1} \scru{w}{gl}_i f(\scru{x}{gl}_i) + \scru{\mathcal{E}}{gl}_n(f) 
\end{equation*}
where 
\begin{equation} \label{eq-approx-error-GL}
\scru{\mathcal{E}}{gl}_n(f) 
= \frac{(n!)^4}{[(2n)!]^3}\cdot \frac{2^{2n+1}}{2n+1}\cdot f^{(2n)}(\eta), \quad \eta\in (-1,1) 
\end{equation}
is the approximation error.
\end{thm}

The standard interval $[-1, 1]$ can be converted to a general finite interval $[a,b]$ by a linear transform $\psi: x \mapsto t = (b-a)x/2 + (b+a)/2$,
which implies that 
\begin{equation} \label{eq-integral-transform-GL}
\int^b_a f(x)\dif x = \frac{b-a}{2}\int^{1}_{-1} f\left(\frac{b+a +(b-a)x}{2} \right) \dif x.
\end{equation}
According to \eqref{eq-compos-deriv}, we can find that
\begin{equation} \label{eq-2n-derivative}
\left[\fracode{}{x}\right]^{2n} f\left(\alpha x + \beta \right)
= \alpha^{2n}f^{(2n)}\left(\alpha x + \beta \right).
\end{equation}
Let $\xi = [b+a + (b-a)\eta]/2$ where $\eta \in (-1,1)$,
then $\xi \in(a,b)$. Substitute \eqref{eq-integral-transform-GL} and \eqref{eq-2n-derivative} into \eqref{eq-approx-error-GL}, we immediately have
\begin{equation}
\begin{split}
\int^b_a f(x)\dif x 
&= \frac{b-a}{2}\sum^n_{i=1} \scru{w}{gl}_i 
 f\left(\frac{b+a+ (b-a)\scru{x}{gl}_i}{2}, k\right) \\
 &  + \left[\frac{b-a}{2}\right]^{2n+1}\cdot\frac{(n!)^4}{[(2n)!]^3}\cdot \frac{2^{2n+1}}{2n+1}f^{(2n)}(\xi)
\end{split}
\end{equation}
from Theorem \ref{thm-approx-error-GL}. 

We now generalize the definite integral $\displaystyle{\int^b_a f(x)\dif x}$ to the abstract function
$\displaystyle{\int^b_a f(x,k)\dif x}$. The deduction above implies the following corollary.

\begin{cor}
For $f(\cdot, k)\in C^{2n}[a,b]$ and fixed parameter $k$, the integral $\displaystyle{\int^b_a f(x,k)\dif x}$ can be computed with the extended Gauss-Legendre formula
\begin{equation}
\begin{split}
&\int^b_a f(x,k)\dif x \\
 &= \frac{b-a}{2}\sum^n_{i=1} \scru{w}{gl}_i 
 f\left(\frac{b+a+ (b-a)\scru{x}{gl}_i}{2}, k\right) \\
 &\quad   + \scru{\mathcal{E}}{gl}_n(f,a,b,k)
\end{split}
\end{equation}
where $\xi\in[a,b]$ and 
\begin{equation}
\scru{\mathcal{E}}{gl}_n(f,a,b,k) 
= \frac{(n!)^4}{[(2n)!]^3}\cdot \frac{(b-a)^{2n+1}}{2n+1}f^{(2n)}(\xi,k)
\end{equation}
is the approximation error.
\end{cor}

Theoretically, if there is an upper bound $M^f_n(k)$ for $f^{(2n)}(\xi,k)$ such that $\abs{f^{(2n)}(\xi,k)}\le M^f_n(k)$ where $\xi \in (a,b)$, then we can find an integer 
\begin{equation}
\begin{split}
\scru{n}{gl}_f 
&= \scru{u}{gl}_f(k,\epsilon,a,b) \\
&=  \arg \min_{n\in \mathbb{N}}\set{n: \frac{(n!)^4}{[(2n)!]^3} \frac{(b-a)^{2n+1}}{2n+1}M^f_n(k) < \epsilon}
\end{split}
\end{equation}
such that $\scru{\mathcal{E}}{gl}_n(f,a,b,k) < \epsilon$ when $n > \scru{n}{gl}_f(k,\epsilon,a,b)$.

\section{Orthogonal Polynomials}
\Tab \ref{tab-orthogonal-polynomials} shows the four typical orthogonal polynomials and the corresponding expressions of $A(n), B(n), C(n), c_0, c_1, c_2$ for the Legendre polynomials $P_n(x)$, Chebyshev polynomials $T_n(x)$, Laguerre polynomials $L_n(x)$ and Hermite polynomials $H_n(x)$ of interest in wide applications. 

\begin{table*}[htp]
\centering
\caption{Parameters for the iterative formula of orthogonal polynomials $Y_n(x)$ such that $Y_0(x) = c_0, Y_1(x) = c_1 x + c_2$ and $Y_n(x) = [A(n)x + B(n)]Y_{n-1}(x) - C(n) Y_{n-2}(x)$ for $n\ge 2$.}
\label{tab-orthogonal-polynomials}
\begin{tabular}{lclccccccc}
\hline
\textbf{Type}  &  $\rho(x)$ & $Y_n(x)$  & \textbf{Domain} &  $A(n) $ & $B(n)$ &  $C(n)$  & $c_0$ & $c_1$  &  $ c_2$ \\
\hline   \hline
Legendre &  $1$ & $P_n(x)=\cfrac{1}{2^nn!}\cdot\cfrac{\dif^n}{\dif x^n}[(x^2-1)^n]$  & $[-1,1]$ & $2-\frac{1}{n}$ & $0$  &  $1-\frac{1}{n}$  & $1$ & $1$  & $0$\\ 
Chebyshev & $\cfrac{1}{\sqrt{1-x^2}}$  & $T_n(x)=\cos(n\arccos x)$ & $[-1,1]$ & $2$              & $0$  & $1$               & $1$ & $2$  & $0$\\ 
Laguerre & $\me^{-x}$ & $L_n(x)=\cfrac{\me^x}{n!}\cdot\cfrac{\dif^n}{\dif x^n}\left[x^n\me^{-x}\right]$ & $[0,+\infty)$ & $-\frac{1}{n}$   & $2-\frac{1}{n}$   & $1-\frac{1}{n}$  & $1$ & $-1$ & $1$\\ 
Hermite  & $\me^{-x^2}$ & $H_n(x)=(-1)^n\me^{x^2}\cfrac{\dif^n}{\dif x^n}\left[\me^{-x^2}\right]$ & $(-\infty, +\infty)$ & $2$              & $0$  & $2(n-1)$    & $1$ & $2$  & $0$\\         
\hline
\end{tabular}
\end{table*}

\section{Nowton's Method for Solving Root}

\Algr~\ref{alg-root-Newton} is used to calculate the root of the non-linear equation $f(x) = 0$ with an initial value $x_0$ via Newton's method (also named by Newton-Raphson method). Essentially, the Newton's method for solving root is a special kind of the fixed-point algorithm in which the contractive mapping is specified by \eqref{eq-Newton-map}. 
Thus what we should do is to define the updating function $\scrd{\mathcal{A}}{Newton}$. Therefore, we first design a procedure \ProcName{NewtonUpdate} for $\scrd{\mathcal{A}}{Newton}$, then we design the procedure \ProcName{RootSolverNewton} by calling the procedure \ProcName{FixedPointSolver}. For the purpose of solving the $n$ roots of Legendre polynomial $P_n(x)$ encapsulated by the procedure \ProcName{CalcRootsPolynomLegendre}, we should set the argument (pointer fo function) \PtrToFun{f} in \ProcName{NewtonUpdate} with the procedure \ProcName{PolynomLegendre} and call the procedure \ProcName{RootSolverNewton} for $n$ times by a loop. \Fig \ref{fig-calc-roots-Legendre-procedures} shows the logic relation of the procedures involved for calculating the $n$ roots of $P_n(x)$. 

\begin{breakablealgorithm}\label{alg-root-Newton}
\caption{Calculate the root of the non-linear equation $f(x) = 0$ with an initial value $x_0$ via Newton's method (\ProcName{RootSolverNewton})} 
\begin{algorithmic}[1]
\Require Pointers to function \PtrToFun{f} for the function $f(x)$, initial value $x_0$, precision $\epsilon$.
\Ensure the root of $f(x) = 0$.
\Function{RootSolverNewton}{\PtrToFun{f}, $x_0$, $\epsilon$}
\State $\PtrToFun{Update}\gets \ProcName{NewtonUpdate} $;
\State $\PtrToFun{Dist}\gets \ProcName{Dist1d} $;
\State $\cpvar{root}\gets \ProcName{FixedPointSolver}(\PtrToFun{Update},
\PtrToFun{Dist},x_0,\epsilon)$;
\State \Return \cpvar{root};
\EndFunction
\end{algorithmic}
\end{breakablealgorithm}

\section{Test Cases for Verification and Validation} 
\label{sec-testcases-VVT}

The optional test cases for verifying and validating the modules involved in constructing the key algorithms and corresponding procedures for computing the CEI-1 are shown in \Tab  \ref{tab-VVT-test-cases}. 

\begin{table*}[htbp]
\centering
\caption{Optional Test Cases for Computational Tasks}
\label{tab-VVT-test-cases}
\begin{tabular}{llll}
\hline
\textbf{Computational Task} & \textbf{Optional Test Cases} & \textbf{Parameters Binding}   &  \textbf{Results for Verification}    \\
\hline \hline
$\displaystyle{S(x) = \sum^\infty_{j=0}c_jx^j}$ 
& $\displaystyle{\cos x = \sum^{\infty}_{j=0}(-1)^j\frac{x^{2j}}{(2j)!}} $
& $\displaystyle{c_{2j} = \frac{(-1)^{j}}{(2j)!}}$     & $\displaystyle{\cos \pi=-1, \dots}$  \\
& $\displaystyle{\sin x = \sum^{\infty}_{j=0}(-1)^j\frac{x^{2j+1}}{(2j+1)!}} $
& $\displaystyle{c_{2j+1} = \frac{(-1)^{j}}{(2j+1)!}}$     & $\displaystyle{\sin \frac{\pi}{2}=1, \cdots}$  \\
& $\displaystyle{\me^{x} = \sum^{\infty}_{j=0}\frac{x^{j}}{j!}} $
& $\displaystyle{c_{j} = \frac{1}{j!}}$     & $\displaystyle{\me^{1}=2.71828\cdots, \cdots}$  \\
& $\displaystyle{\arctan x = \sum^{\infty}_{j=0}(-1)^j\frac{x^{2j+1}}{2j+1}} $
& $\displaystyle{c_{2j+1} = \frac{(-1)^j}{2j+1}}$     & $\displaystyle{\arctan 1 = \frac{\pi}{4}}$  \\
\hline
$P_n(x)$  & $P_0(x)$  & $n = 0$      &   $P_1(x) = 1 $ \\
          & $P_1(x)$  & $n = 1$      &   $P_1(x) = x $ \\
          & $P_2(x)$  & $n = 2$      &   $\displaystyle{P_2(x) = (3x^2-1)/2}$  \\
          & $P_3(x)$  & $n=3$        &   $\displaystyle{P_3(x) = (5x^3-3x)/2}$\\  
\hline
$P'_n(x)$  & $P'_0(x)$  & $n = 0$     &   $P'_1(x) = 0 $ \\
          & $P'_1(x)$  & $n = 1$      &   $P'_1(x) = 1 $ \\
          & $P'_2(x)$  & $n = 2$      &   $\displaystyle{P'_2(x) = 3x}$  \\
          & $P'_3(x)$  & $n=3$        &   $\displaystyle{P'_3(x) = (15x^2-3)/2}$\\  
\hline
$f(x,\alpha)=0$ & $\displaystyle{f(x) =\prod^{n}_{i=1}(x-\gamma_i)=0}$ & $\gamma_1<\gamma_2<\cdots<\gamma_n$ & $x \in \set{\gamma_1, \cdots, \gamma_n}$\\
$P_n(x)=0$ & $P_1(x)=0$  & $n = 1$      &   $x = 0 $ \\
          & $P_2(x)=0$   & $n = 2$      &   $x\in  \set{-\frac{\sqrt{3}}{3}, \frac{\sqrt{3}}{3}}$  \\
          & $P_3(x)=0$  & $n=3$        &   $x\in \set{-\sqrt{\frac{3}{5}}, 0, \sqrt{\frac{3}{5}}}$\\  
\hline
$\displaystyle{\scrd{F}{C}(k)=\int^a_{-a}\frac{f(x,k)}{\sqrt{a^2-x^2}}\dif x}$ &  $\displaystyle{\int^a_{-a}\frac{1}{\sqrt{1-x^2}}\dif x}$    &$f(x,k) = 1$  
& $\displaystyle{ 2\arcsin a} $ \\
 &  $\displaystyle{\int^1_{-1}\frac{\cos kx}{\sqrt{1-x^2}}\dif x}$    &$f(x) = \cos kx$  
& $\displaystyle{J_0(k)}$ (Bessel function) \\
 &  $\displaystyle{\int^1_{-1}\frac{\me^{kx}}{\sqrt{1-x^2}}\dif x}$    &$f(x) = \me^{kx}$  
& $\displaystyle{\pi\me^{kx}} $ \\
\hline
$\displaystyle{\scrd{F}{L}(k)=\int^b_af(x,k)\dif x}$ & $\displaystyle{\int^b_{a}\cos(x+k)\dif x}$ & $ f(x,k) = \cos(x+k)$   & $ \sin(k+b) - \sin (k+a)$ \\
   &  $\displaystyle{\int^b_a\frac{(k+x)^2}{3}\dif x}$  &  $\displaystyle{f(x,k)=\frac{(x+k)^2}{3}}$    & $(k+b)^3-(k+a)^3$   \\  
 &$\displaystyle{\int^b_a \me^{kx}\dif x }$  &$\displaystyle{f(x,k)=\me^{kx}}$    & $\displaystyle{\frac{\me^{kb} - \me^{ka}}{k}}$   \\  
 \hline
$\AGM(a,b)$   &  $\displaystyle{\AGM(2,3)}$   & $a=2,b=3$   & $\displaystyle{\frac{\pi}{2}/\int^{\frac{\pi}{2}}_0\frac{\dif \theta}{\sqrt{2^2\cos^2\theta + 3^2\sin^2\theta}}}$  \\
              &  $\displaystyle{\frac{1}{\AGM(1+x,1-x)}}$ & $x\in[0,1)$ & $\displaystyle{\sum^\infty_{j=0}\left[\frac{(2j-1)!!}{(2j)!!}\right]^2x^{2j}}$ \\
              &  $\AGM(a,b)$   &  $a,b\in \Real^+$  & $\AGM(b,a)=\AGM\left(\cfrac{a+b}{2}, \sqrt{ab}\right)$   \\
\hline
\end{tabular}
\end{table*}

\Tab \ref{tab-VVT-test-cases} lists some optional test cases for the computational tasks in computing the CEI-1. We give some interpretations here:
\begin{itemize}
\item For computing the infinite series $S(x) = \sum^\infty_{j=0} c_m x^m$, we can use some typical functions (say $\cos x, \sin x, \me^x$ and $\arctan x$)  and their expansions of series to verify the numerical solutions obtained by \Algr \ref{alg-InfinteSeriesSolver} by comparing their analytical solutions.  
\item For generating the Legendre polynomial $P_n(x)$ and its derivative $P'_n(x)$ , we can use the expressions of $P_0(x), P_1(x), P_2(x)$ and $P_3(x)$ to test whether \Algr \ref{alg-OrthoPolynom}, \Algr \ref{alg-DerivOrthoPolynom}, \Algr \ref{alg-PolynomLegendre}, \Algr \ref{alg-DerivPolynomLegendre} and the corresponding small procedures involved work well. 
\item For computing the root of nonlinear equation $f(x,\alpha)=0$ with \Algr \ref{alg-AbstractFixedPointSolver}, \Algr \ref{alg-FixedPointSolver} and \Algr \ref{alg-NewtonUpdate}, we can construct 
$\displaystyle f(x) = \prod^n_{i=1}(x-\gamma_i)$
 to test the numerical solution to the analytical solution $x_i = \gamma_i$ for positive $n$ and $\gamma_1<\gamma_2< \cdots < \gamma_n$, say $\gamma_i = i$.  
\item For calculating the $n$ roots of $P_n(x)$ with Newton's method, we can use the analytical solution to $P_1(x) = x = 0, P_2(x) =(3x^2-1)/2= 0, P_3(x)=x(5x^2-3)/2=0$ to verify the numerical solution obtained by \Algr \ref{alg-RootSolverOrthPolynom}, \Algr \ref{alg-NewtonUpdate} and \Algr \ref{alg-CalcRootsLegendre}.  
\item For computing the $\displaystyle{\scrd{F}{C}(k)=\int^a_{-a}\frac{f(x,k)}{\sqrt{a^2-k^2}}\dif x}$ with \Algr \ref{alg-IntegralSolverGC}, we can put $f(x,k) \in \set{1, \cos kx, \me^{kx}}$ and compare the numerical solutions with their analytical counterparts
\begin{equation}
\left\{
\begin{split}
\int^1_{-1}\frac{1}{\sqrt{1-x^2}}\dif x &= \pi \\
\int^1_{-1}\frac{\cos kx}{\sqrt{1-x^2}}\dif x &= J_0(k)\\
\int^1_{-1}\frac{\me^{kx}}{\sqrt{1-x^2}}\dif x &= \pi\me^{kx}
\end{split}
\right.
\end{equation} 
where $J_0(\cdot)$ is the $0$-th order Bessel function of the first kind.
\item For solving the $\displaystyle{\scrd{F}{L}(k) = \int^b_a f(x,k)\dif x}$  with \Algr \ref{alg-IntegralSolverGL}, we can set $$f(x,k)\in \set{\cos (x+k), \frac{1}{3}(x+k)^2, \me^{kx}}$$
to check the numerical solutions. Note that 
\begin{align*}
\lim_{k\to 0} \int^b_a \me^{kx}\dif x 
&= \lim_{k\to 0} \frac{\me^{kb}-\me^{ka}}{k} \\
&= b-a = \int^b_a \lim_{k\to 0} \me^{kx} \dif x
\end{align*}
by the Newton-Leibnitz formula and the rule of L'Hospital in calculus. This result also can be verified with numerical integration.  
\item The verification of \Algr \ref{alg-AGM} for solving $\AGM(a,b)$ is challenging if compared with the other computational tasks since there is no simple analytical expression for it. However, as illustrated at the bottom of \Tab \ref{tab-VVT-test-cases}, we can check the value of $\AGM(a,b)$ through cross-verification by the symmetric-iterative property and the relation with infinite series and integral according to equations \eqref{eq-AGM-symmetric},  \eqref{eq-AGM-iterative}, \eqref{eq-AGM-series} and \eqref{eq-AGM-integral} since the integral and infinite series can be verified with available methods as mentioned above. 
\end{itemize}

\section*{Declaration of competing interest}

The authors declare that they have no known competing financial interests or personal relationships that could have appeared to influence the work reported in this paper.

\section*{Code \& data availability statement}
The algorithms developed in this paper can be obtained at the GitHub site 
\url{https://github.com/GrAbsRD/CEI-1/tree/main} 
where a long table with filename \lstinline|Kvalues.txt| for the values of $K(k)$ is provided and the reader can run the code to verify the algorithms and reuse the code according to his/her interests.

\section*{Acknowledgment}

This work was supported partly by the National Natural Science Foundation of China under grant number 62167003, and partly by the Hainan Provincial Education and Teaching Reform Project of Colleges and Universities under grant number Hnjg2021-37.


\bibliographystyle{unsrt}

\begin{thebibliography}{10}

\bibitem{TableISP2015}
I.~S. Gradshteyn and I.~M. Ryzhik.
\newblock {\em Table of Integrals, Series, and Products}.
\newblock Academic Press (Elsevier Inc), London, eighth edition, 2015.
\newblock Translated from Russian by Scripta Technica, Inc.

\bibitem{Abramo1965}
M.~Abramowitz and I.A. Stegun.
\newblock {\em Handbook of Mathematical Functions: with Formulas, Graphs, and Mathematical Tables}.
\newblock Dover Publications, New York, 2nd edition, 1965.
\newblock chapter 17.

\bibitem{HYZhang2024CdioCt}
Hong-Yan Zhang, Yu~Zhou, Yu-Tao Li, Fu-Yun Li, and Yong-Hui Jiang.
\newblock {CDIO-CT} collaborative strategy for solving complex {STEM} problems
  in system modeling and simulation: {A}n illustration of solving the period of
  mathematical pendulum.
\newblock {\em Computer Applications in Engineering Education}, 32(2):e22698,
  2024.
\newblock \url{https://doi.org/10.1002/cae.22698}.

\bibitem{Atkinson2009}
Kendall Atkinson and Wei-Min Han.
\newblock {\em Theoretical Numerical Analysis: A Functional Analysis
  Framework}, volume~39 of {\em Texts in Applied Mathematics}.
\newblock Springer, New York, 3rd edition, 2009.

\bibitem{TAOCP1}
Donald~E. Knuth.
\newblock {\em The Art of Computer Programming}, volume 1: Fundamental
  Algorithms.
\newblock Addison-Wesley Professional, New York, 3rd edition, 1997.

\bibitem{King1924AGM}
L.~V. King.
\newblock {\em On the Direct Numerical Calculation of Elliptic functions and
  Integrals}.
\newblock Cambridge University Press, Cambridge, 1924.

\bibitem{Borwein1987}
Jonathan~M. Borwein and Peter~B. Borwein.
\newblock {\em {Pi and the AGM: A Study in Analytic Number Theory and
  Computational Complexity}}.
\newblock Wiley-Interscience and Canadian Mathematics Series of Monographs and
  Texts. Wiley-Interscience, 1987.

\bibitem{CourantHilbert-vol1}
Richard Courant and David Hilbert.
\newblock {\em Methods of Mathematical Physics}, volume~1.
\newblock Wiley-VCH, New York, 2nd edition, 1989.

\bibitem{SJZhang1996}
Shan-Jie Zhang and Jian-Ming Jin.
\newblock {\em Computation of Special Functions}.
\newblock Wiley-Interscience, New York, 1996.
\newblock the Chinese version was published in 2011 by the Nanjing University
  Press in China.

\bibitem{HandbookMPES2011}
A.~I. Chernoutsan, A.~V. Egorov, A.~V. Manzhirov, A.~D. Polyanin, V.~D.
  Polyanin, V.~A. Popov, B.~V. Putyatin, Yu.~V. Repina, V.~M. Safrai, and A.I.
  Zhurov.
\newblock {\em A Concise Handbook of Mathematics, Physics, and Engineering
  Sciences}.
\newblock CRC Press, Boca Roton, 2011.

\bibitem{Brooks2010}
Frederick~P. Brooks.
\newblock {\em {The Design of Design: Essays from A Computer Scientist}}.
\newblock Addison-Wesley Professional, 2010.

\end{thebibliography}

\noindent\textbf{Citation}:
Hong-Yan Zhang, Wen-Juan Jiang, Open source implementations of numerical algorithms for computing the complete elliptic integral of the first kind, \textit{Results in Applied Mathematics}, 2024, 23(8): e100479, \url{https://doi.org/10.1016/j.rinam.2024.100479}.

\end{document}